\documentclass{ws-m3as-ppn}

\usepackage[T1]{fontenc}
\usepackage[utf8]{inputenc}
\usepackage{amsmath,amsfonts,amssymb,amsxtra,bbm}
\usepackage{fourier,setspace,graphicx,color,pdflscape}
\usepackage[colorlinks=true]{hyperref}
\hypersetup{urlcolor=blue, linkcolor=blue, citecolor=red, anchorcolor=blue]}
\usepackage{silence}

\newcommand{\R}{{\mathbb R}}
\renewcommand{\S}{{\mathbb S}}
\newcommand{\N}{{\mathbb N}}

\newcommand{\be}[1]{\begin{equation}\label{#1}}

\newcommand{\ee}{\end{equation}}
\renewcommand{\(}{\left(}
\renewcommand{\)}{\right)}
\newcommand{\ird}[1]{\int_{\R^d}{#1}\,dx}

\newcommand{\isd}[1]{\int_{\S^d}{#1}\,d\sigma}
\newcommand{\nrm}[2]{\left\|{#1}\right\|_{#2}}
\newcommand{\nrmS}[2]{\|{#1}\|_{\mathrm L^{#2}(\S^d)}}
\newcommand{\nrmR}[2]{\|{#1}\|_{\mathrm L^{#2}(\R^d)}}
\newcommand{\seq}[2]{({#1}_{#2})_{#2\in\N}}

\newcommand{\nrms}[2]{\|{#1}\|_{\mathrm L^{#2}(\S^1)}}
\newcommand{\is}[1]{\int_0^{2\,\pi}{#1}\,d\sigma}
\newcommand{\iE}[1]{\int_{I(E)}{#1}}

\makeatletter
\@namedef{subjclassname@2020}{\textup{2020} Mathematics Subject Classification}
\makeatother
\newcommand{\msc}[1]{\href{https://zbmath.org/classification/?q=cc:#1}{#1}}

\definecolor{darkgreen}{rgb}{.0,.4,.2}

\newcommand{\Pot}{\mathcal V}

\usepackage{etoolbox}
\newcounter{taggedeq}
0
\pretocmd{\equation}{\stepcounter{taggedeq}}{}{}

\begin{document}
\title{Interpolation inequalities on the sphere and phase transition:\\ rigidity, symmetry and symmetry breaking}

\author{Esther Bou Dagher}
\address{CEREMADE (CNRS UMR n$^\circ$ 7534),\\ PSL university, Universit\'e Paris-Dauphine,\\ Place de Lattre de Tassigny, 75775 Paris 16, France.\\
E-mail: \email{esther.bou-dagher@dauphine.psl.eu}}

\author{Jean Dolbeault}
\address{CEREMADE (CNRS UMR n$^\circ$ 7534),\\ PSL university, Universit\'e Paris-Dauphine,\\ Place de Lattre de Tassigny, 75775 Paris 16, France.\\
E-mail: \email{dolbeaul@ceremade.dauphine.fr}}

\maketitle
\thispagestyle{empty}

\begin{abstract}
This paper is devoted to the study of phase transitions associated to a large family of Gagliardo-Nirenberg-Sobolev interpolation inequalities on the sphere depending on one parameter. We characterize symmetry and symmetry breaking regimes, with a phase transition that can be of first or second order. We establish various new results and study the qualitative properties of the branches of solutions to the Euler-Lagrange equations.
\end{abstract}

\keywords{Interpolation; Gagliardo-Nirenberg-Sobolev inequalities; best constants; branches of solutions; optimal functions; bifurcations; symmetry; symmetry breaking; uniqueness; rigidity; phase transition.}

\ccode{AMS Subject Classification: \msc{58J05},\msc{35B06}, \msc{26D10}.}
\section{Introduction}\label{Sec:Intro}

On the $d$-dimensional unit sphere, the classical \emph{Gagliardo-Nirenberg-Sobolev interpolation inequality},
\be{GNS}\tag{GNS}
\nrmS{\nabla u}2^2\ge\frac d{p-2}\(\nrmS up^2-\nrmS u2^2\)
\ee
for any $u\in\mathrm H^1(\S^d,d\sigma)$, was proved with optimal constant by M.-F.~Bidaut-V\'eron and L.~V\'eron~\cite{1991}, and also by W.~Beckner~\cite{MR1230930} using another method, for $p>2$. The case $p<2$ and some cases with $p>2$ are covered by the \emph{carr\'e du champ} method coupled with the heat flow~\cite{Bakry-Emery85,MR808640,1504}. Here $d\ge1$ is an integer, $d\sigma=\big|\mathbb S^d\big|^{-1}\,dv_g$ denotes the uniform probability measure on the unit sphere $\S^d$, and $dv_g$ can be seen as the measure induced on $\S^d\subset\R^{d+1}$ by Lebesgue's measure on $\R^{d+1}$. The space $\mathrm L^q(\S^d,d\sigma)$ with $q\in[1,+\infty)$ is the standard Lebesgue space with norm
\[
\nrmS uq:=\(\isd{|u|^q}\)^{1/q}.
\]
The fact that $d\sigma$ is a probability measure means that $\nrmS uq=1$ if $u=1$ a.e.~on $\S^d$ and we also have $\nrmS u{q_1}\le\nrmS u{q_2}$ for any $u\in\mathrm L^{q_2}(\S^d,d\sigma)$ as soon as $1\le q_1\le q_2$. The space $\mathrm H^1(\S^d,d\sigma)$ is obtained by completion of the space of smooth functions on $\S^d$ with respect to the norm $u\mapsto\big(\nrmS{\nabla u}2^2+\nrmS u2^2\big)^{1/2}$. Inequality~\eqref{GNS} holds in the range
\[
p\in[1,2)\cup(2,\infty)\quad\mbox{if}\quad d=1,\,2\quad\mbox{and}\quad p\in[1,2)\cup(2,2^*]\quad\mbox{if}\quad d\ge3\,,\quad\mbox{with}\quad 2^*:=\frac{2\,d}{d-2}\,.
\]
By convention, we take $2^*=+\infty$ if $d=1$, $2$. The \emph{carr\'e du champ} method has been extended to nonlinear diffusion flows~\cite{MR2381156,Dolbeault20141338} in order to cover the whole range of admissible exponents on manifolds with positive curvature, with applications to the sphere~\cite{DEKL2012,DEKL2014,DolEsLa-APDE2014}. Inequality~\eqref{GNS} can be considered as a special case of the family of interpolation inequalities in $\mathrm H^1(\S^d,d\sigma)$
\be{Ineq:GNS1}
\nrmS{\nabla u}2^2+\frac\lambda{p-2}\,\nrmS u2^2\ge\frac{\overline{\mu}(p,\lambda)}{p-2}\,\nrmS up^2
\ee
which depends on a positive parameter $\lambda$, with optimal constant $\overline{\mu}(p,\lambda)$. One can indeed check that
\[
d=\max\left\{\lambda>0\,:\,\overline{\mu}(p,\lambda)=\lambda\right\}\,.
\]
By construction, the function $\lambda\mapsto\overline{\mu}(p,\lambda)$ is monotone non-decreasing and concave. Branches of critical functions, parametrized by $\lambda$, which realize the equality case in~\eqref{Ineq:GNS1}, have been characterized by the \emph{carr\'e du champ} method~\cite{DolEsLa-APDE2014,1504,Dolbeault_2020b}. In this paper, we investigate the larger family of inequalities
\be{Ineq:GNS2}
\(\nrmS{\nabla u}2^2+\frac\lambda{p-2}\,\nrmS u2^2\)^\theta\Bigg(\nrmS u2^2\Bigg)^{1-\theta}\ge\(\frac{\mu(p,\theta,\lambda)}{p-2}\)^\theta\,\nrmS up^2\,,
\ee
for any $u\in\mathrm H^1(\S^d,d\sigma)$, where $\theta>0$ and $\lambda>0$ are two parameters and $\mu(p,\theta,\lambda)$ is the corresponding optimal constant. Depending on the value of $\theta$, we observe several scenarios of \emph{phase transition}. The case $\theta=1$ corresponds to~\eqref{Ineq:GNS1} and the range $\theta>1$ is admissible. The range $\theta<1$ deserves some additional care, as there is another limitation, $\theta\ge\theta_\star$, where
\[
\theta_\star:=d\,\frac{p-2}{2\,p}
\]
is the exponent in the \emph{Euclidean Gagliardo-Nirenberg-Sobolev inequality}
\be{Ineq:Rd}
\nrmR{\nabla f}2^{\theta_\star}\,\nrmR f2^{1-\theta_\star}\ge\mathcal C_{\mathrm{GNS}}(p)\,\nrmR fp\quad\forall\,f\in\mathrm H^1(\R^d,dx)\,.
\ee
Our purpose is to study $\lambda\mapsto\mu(p,\theta,\lambda)$ and characterize whenever possible the optimal functions. Let us define the \emph{Bakry-Emery exponent} $2^\#$, the parameters $\gamma^{\#}$ and $\theta^\#$, such that $2<2^\#\le2^*$ and $\theta_\star<\theta^\#<1$, by
\begin{align*}
&2^\#:=\frac{2\,d^2+1}{(d-1)^2}\,,\quad\gamma^{\#}:=\(\frac{d-1}{d+2}\)^2(p-1)\,(2^\#-p)\quad\mbox{and}\quad\,\theta^\#:=\(\frac{\gamma^{\#}}{p-2}+1\)^{-1}&\quad\mbox{if}\quad d\ge2\,,\\
&2^\#=+\infty\,,\quad\gamma^{\#}=\frac{p-1}3\quad\mbox{and}\quad\,\theta^\#:=3\,\frac{p-2}{4\,p-7}&\quad\mbox{if}\quad d=1\,.
\end{align*}
Inequality~\eqref{Ineq:GNS2} is known, for any $p\in(2,2^\#)$ but only for $\theta=\theta^\#$, from Inequality~(2.4) of~\cite{Dolbeault_2020b}. This inequality arises as a consequence of improvements of~\eqref{Ineq:GNS1} based on the \emph{carr\'e du champ} method, written in the form
\be{Ineq:GNSimproved}
\nrmS{\nabla u}2^2\ge\frac{d\,\theta^\#}{p-2}\(\nrmS up^{2/\theta^\#}\,\nrmS u2^{2\,(\theta^\#-1)/\theta^\#}-\nrmS u2^2\)
\ee
for any $u\in\mathrm H^1(\S^d,d\sigma)$. Our first result deals with Inequality~\eqref{Ineq:GNS2} and properties of $\mu(p,\theta,\lambda)$.\vspace*{-8pt}
\begin{theorem}\label{Thm1} Let $d\ge1$, $p\in(2,+\infty)$ if $d=1$, $2$, $p\in(2,2^*]$ if $d\ge3$, and $\theta\ge\theta_\star$. Inequality~\eqref{Ineq:GNS2} holds for some optimal $\mu(p,\theta,\lambda)>0$, for any $\lambda>0$ and equality is achieved if $p<2^*$ and $\theta>\theta_\star$. The function $\lambda\mapsto\mu(p,\theta,\lambda)$ is monotone non-decreasing, monotone increasing if $\theta>\theta_\star$, concave and such that, for any $\lambda>0$,
\begin{align*}
&\mu(p,\theta,\lambda)\le\lambda\quad\mbox{and}\quad\mu(p,\theta,\lambda)<\lambda\quad\mbox{if}\quad\lambda>d\,\theta\,,\\
&\mu(p,\theta,\lambda)=\lambda\quad\mbox{if}\quad\lambda\le d\,\theta\,,\quad\theta\ge\theta^\#\quad\mbox{and}\quad p\in(2,2^\#]\quad\mbox{if}\quad d\ge2\quad\mbox{or}\quad p>2\quad\mbox{if}\quad d=1\,.
\end{align*}
Moreover, for any $\theta\ge\theta_\star$, there is some explicit $\kappa>0$ such that
\be{asymptotics}
\mu(p,\theta,\lambda)\sim\kappa\,\lambda^{1-\theta_\star/\theta}\quad\mbox{as}\quad\lambda\to+\infty\,.
\ee
\end{theorem}
\noindent Whenever $\mu(p,\theta,\lambda)=\lambda$, constant functions realize the equality case in~\eqref{Ineq:GNS2}. If these constant functions are the only optimal functions, we shall say that there is \emph{symmetry} because none of the optimal functions points in a privileged direction. The equality case is achieved by non-constant functions if $\mu(p,\theta,\lambda)<\lambda$, in which case rotations and multiplications by a constant generate a $d$-dimensional set of optimal functions. In that case, we shall say that there is \emph{symmetry breaking}. The key issue is to understand for which value of~$\lambda$ the \emph{phase transition} from a regime with symmetry to a regime of symmetry breaking occurs, that is, to determine
\[
\lambda_{\mathrm T}(p,\theta):=\inf\big\{\lambda>0\,:\,\mu(p,\theta,\lambda)<\lambda\big\}\,.
\]
We read from Theorem~\ref{Thm1} that $\lambda_{\mathrm T}(p,\theta)\le d\,\theta$, and $\lambda_{\mathrm T}(p,\theta)=d\,\theta$ if $\theta\ge\theta^\#$ and $p\in(2,2^\#)$. Our second statement goes as follows. Explicit estimates of $\lambda_{\mathrm T}$ will be given later, but the general picture goes as follows.
\begin{theorem}\label{Thm2} Let $d\ge1$, $p\in(2,+\infty)$ if $d=1$, $2$, $p\in(2,2^*]$ if $d\ge3$ and $\theta\ge\theta_\star$. Then $\lambda_{\mathrm T}(p,\theta)$ is positive and $\lambda_{\mathrm T}(p,\theta)<d\,\theta$ if $\theta-\theta_\star\ge0$ is taken small enough and $p<p_\ast$ for some $p_\ast>2$.\end{theorem}

\medskip Inequalities~\eqref{Ineq:GNS2} provide us with a simple model of nonlinear \emph{phase transition}, which fits into the \emph{Ehrenfest classification}~\cite{Jaeger_1998}. For given values of $p$ and $\theta$, the transition from a regime with symmetry to a regime of symmetry breaking is driven by the parameter $\lambda$. The curve $\lambda\mapsto\mu(p,\theta,\lambda)$ can either be smooth (\emph{second order phase transition}), or give rise to a jump in the derivative with respect to $\lambda$ (\emph{first order phase transition}). It turns out that a detailed description can be obtained in the simplest case, $d=1$, which is particularly enlightening: see Section~\ref{Sec:d=1} for details and also Section~\ref{Sec:Numerics} for various numerical illustrations. As we shall see with less details, similar patterns are also present in higher dimensions.

A key tool for the analysis of~\eqref{Ineq:GNS1} is a generalized \emph{carr\'e du champ} method which amounts to evolve the \emph{deficit} in~\eqref{Ineq:GNS1}, that is, the difference of the two sides of the inequality, by a nonlinear diffusion equation. This allows us to reduce the issue of the optimality from a global nonlinear variational problem to a local problem that can be studied with spectral tools. Such a strategy enters in the class of the \emph{entropy methods}. From the point of view of physics, the evolution equation that we use is a diffusion at fixed temperature and the notion of (generalized) entropy is in fact a \emph{free energy}. 

The \emph{symmetry versus symmetry breaking} issue in~\eqref{Ineq:GNS2} provides us with a rather simple mathematical model, in which both types of phase transition occur depending on the choice of the parameter $\theta$. Similar results on phase transitions have been observed in \emph{Caffarelli-Kohn-Nirenberg inequalities}~\cite{Caffarelli1984,Catrina2001}, a well known example of functional inequalities in which symmetry breaking plays an important role, but which also raises various technical difficulties and is for this reason less understood~\cite{DEL-JEPE,BDNS2023}. The sharp threshold for the symmetry range of the parameters in such inequalities is known in the critical case~\cite{Felli2003,MR3570296} and in the subcritical case~\cite{Bonforte2017a,Dolbeault2017}, under a condition which corresponds to $\theta=1$ in~\eqref{Ineq:GNS2}. However, a complete parabolic proof based on entropy methods is so far missing and results rely on \emph{rigidity} results: see~\cite{DEL-JEPE,Dolbeault_2022,BDNS2023} for partial results based on flows. The analogy between our interpolation inequalities on the sphere and Caffarelli-Kohn-Nirenberg inequalities can be pushed further, as one can study branches of critical functions as in~\eqref{Ineq:GNS1}, see for instance~\cite{Catrina2001,0902,1005}. There is also a family of inequalities involving a parameter $\theta<1$ as in~\eqref{Ineq:GNS2} which has several qualitative properties~\cite{delPino20102045,springerlink:10.1007/s00526-011-0394-y,FreefemDolbeaultEsteban,0951-7715-27-3-435,DoEsFiTe2014} that are similar to the properties that we can observe on the sphere. However Inequality~\eqref{Ineq:GNS2} provides us with a much more detailed example of the phenomenon of phase transition, especially in dimension $d=1$, in a simple functional framework. Moreover, the regularity of the solution of the nonlinear flow on $\S^d$ raises no difficulty~\cite{MR2459454} (to be compared with~\cite{DEL-JEPE,Dolbeault_2022,BDNS2023}).

Let us quote some additional entry points in the literature. M.-F.~Bidaut-V\'eron and L.~V\'eron proved~\cite{1991} Inequality~\eqref{Ineq:GNS1} as a consequence of a \emph{rigidity} result, \emph{i.e.}, the property that there is no other positive critical point than the constant functions. This strategy goes back to the work of B.~Gidas and J.~Spruck~\cite{MR615628}. The proof of~\eqref{Ineq:GNS1} by W.~Beckner~\cite{MR1230930} relies on spectral methods, the Funk-Hecke formula and a duality result introduced by E.~Lieb~\cite{Lieb-83}. An earlier version corresponding to the range $p\in(2,2^\#)$ was established by D.~Bakry and M.~Emery~\cite{Bakry-Emery85,MR808640}, using the \emph{carr\'e du champ} method and the heat flow. We refer to the book of D.~Bakry, I.~Gentil and M.~Ledoux~\cite{MR3155209} for a general overview of the \emph{carr\'e du champ} method in the context of Markov processes or linear diffusion equations, and to several papers~\cite{MR2381156,Dolbeault20141338,1504} for the extension to nonlinear diffusion equations. Inequality~\eqref{Ineq:GNS2} is less standard than~\eqref{Ineq:GNS1} and we are not aware of any specific reference, except for $\theta=\theta^\#$, as a consequence of improvements of~\eqref{Ineq:GNS1} based on the \emph{carr\'e du champ} method~\cite{Dolbeault_2020b}. We also refer to a presentation~\cite{Dolbeault_2020b} of \emph{entropy methods} based on nonlinear flows and \emph{improved interpolation inequalities}, with further references therein. An important motivation for~\eqref{Ineq:GNS2} arises from pure states in Lieb-Thirring estimates and interpolation inequalities for systems~\cite{BDT}. We finally refer to review papers~\cite{1703,Dolbeault:2021wb} on, respectively, functional inequalities and branches of solutions in various frameworks, and entropy methods associated with nonlinear diffusion equations.

\medskip This paper is organized as follows. In Section~\ref{Sec:Prelim}, we state various preliminary results, give some explicit lower bounds on $\mu(p,\theta,\lambda)$ and obtain estimates for the symmetry range. In Section~\ref{Sec:More}, we establish the asymptotics of $\mu(p,\theta,\lambda)$ as $\lambda\to+\infty$ and reparametrize the set of solutions with $\theta\neq1$ by the solutions of the problem with $\theta=1$. How solutions behave near the bifurcation point is studied at formal level in Section~\ref{Sec:Taylor}. Section~\ref{Sec:d=1} is devoted to the case of dimension $d=1$. All non-constant solutions of the Euler-Lagrange equations can be parametrized by an energy, which provides us with various estimates. In Section~\ref{Sec:Numerics}, we collect numerical results which illustrate our main results and explain why phase transitions are of particular interest in the study of Inequality~\eqref{Ineq:GNS2}.

\section{Preliminary results}\label{Sec:Prelim}

Let us consider the optimal constant for~\eqref{Ineq:GNS2}, which is given by the variational problem
\be{mu-Q}
\mu(p,\theta,\lambda)=\inf_{u\in\mathrm H^1(\S^d,d\sigma)\setminus\{0\}}\mathcal Q_{p,\theta,\lambda}[u]\quad\mbox{with}\quad\mathcal Q_{p,\theta,\lambda}[u]:=\frac{(p-2)\,\nrmS{\nabla u}2^2+\lambda\,\nrmS u2^2}{\nrmS up^{2/\theta}\,\nrmS u2^{2-2/\theta}}\,.
\ee

\subsection{Some simple estimates}\label{Sec:Simple}

By using $u\equiv1$ as test function in~\eqref{mu-Q}, we have an easy upper bound on $\mu(p,\theta,\lambda)$.
\begin{lemma}\label{Lem:mu} Let $d\ge1$, $p\in(2,+\infty)$ if $d=1$, $2$, $p\in(2,2^*]$ if $d\ge3$ and $\theta\ge\theta_\star$. For all $\lambda>0$, we have
\[
\mu(p,\theta,\lambda)\le\lambda\,.
\]
\end{lemma}
\noindent Lower bounds are not as simple.
\begin{lemma}\label{Lem:theta} Let $d\ge1$, $p\in(2,+\infty)$ if $d=1$, $2$, $p\in(2,2^*]$ if $d\ge3$. Assume that $\theta_\star\le\theta_1<\theta_2$. Then we have
\be{Mon1}
\mu(p,\theta_\star,\lambda)\le\mu(p,\theta_1,\lambda)\le\mu(p,\theta_2,\lambda)\,.
\ee
Additionally, if $\theta_\star\le\theta_1<\theta_2<p\,\theta_1/(p-2)$ and $\theta_2\le\theta_1/\theta_\star$ if $d\ge3$, then we have
\be{Mon2}
\mu(p,\theta_1,\lambda)\ge\frac1{2\, {\theta_2}}\,\big(p\,\theta_1-(p-2)\,\theta_2\big)\,\mu\(\frac{2\,p\,\theta_1}{p\,\theta_1-(p-2)\,\theta_2},\theta_2,\frac{2\,\lambda\,\theta_2}{p\,\theta_1-(p-2)\,\theta_2}\)\,.
\ee
\end{lemma}
\begin{proof} Since $d\sigma$ is a probability measure, we have $\nrmS u2^2\le\nrmS up^2$ by H\"older's inequality so that
\[
\mathcal Q_{p,\theta_1,\lambda}[u]=\((p-2)\,\frac{\nrmS{\nabla u}2^2}{\nrmS u2^2}+\lambda\)\(\frac{\nrmS u2^2}{\nrmS up^2}\)^{1/\theta_1}\le\((p-2)\,\frac{\nrmS{\nabla u}2^2}{\nrmS u2^2}+\lambda\)\(\frac{\nrmS u2^2}{\nrmS up^2}\)^{1/\theta_2}=\mathcal Q_{p,\theta_2,\lambda}[u]
\]
which proves~\eqref{Mon1}. Concerning~\eqref{Mon1}, we may use H\"older's inequality $\nrmS up\le\nrmS uq^\theta\,\nrmS u2^{1-\theta}$, that~is,
\be{Holder:p-q}
\frac{\nrmS u2^2}{\nrmS up^2}\ge\(\frac{\nrmS u2^2}{\nrmS uq^2}\)^\theta\quad\mbox{where}\quad\theta=\frac{q\,(p-2)}{p\,(q-2)}\,.
\ee
If $\theta>(p-2)/p$, this determines $q=q(p,\theta)$ given by
\be{qptheta}
q(p,\theta):=\frac{2\,p\,\theta}{2-p\,(1-\theta)}\,.
\ee
in the range $(p,2^*)$ with $q=2^*$ in the limit case $\theta=\theta_\star=d\,(p-2)/(2\,p)$, whenever $d\ge3$, and $q$ in the range $(p,\infty)$ if $d=1$ or $d=2$. Let
\be{Lambda}
\Lambda(p,\theta,\lambda):=\frac{q(p,\theta)-2}{p-2}\,\lambda=\frac{2\,\lambda}{2-p\,(1-\theta)}\,.
\ee
Altogether, with $q(p,\theta)=q(p,\theta)$ and $\Lambda=\Lambda(p,\theta,\lambda)$, we can then write that
\[
\mathcal Q_{p,\theta_1,\lambda}[u]\ge\frac{p-2}{q-2}\((q-2)\,\frac{\nrmS{\nabla u}2^2}{\nrmS u2^2}+\Lambda\)\(\frac{\nrmS u2^2}{\nrmS uq^2}\)^{\theta/\theta_1}=\(1-\frac p2\,(1-\theta)\)\mathcal Q_{q,\theta_1/\theta,\Lambda}[u]\,.
\]
If $d\ge3$, the condition $q(p,\theta)\le2\,d/(d-2)$ amounts to $\theta\ge\theta_\star$, which proves~\eqref{Mon2} using $\theta=\theta_1/\theta_2$.
\end{proof}

\begin{corollary}\label{Cor:mu} Let $d\ge1$ and $p\in(2,+\infty)$ if $d=1$, $2$, $p\in(2,2^*]$ if $d\ge3$.\\
If $\theta\ge1$, then
\be{C1}
\overline{\mu}(p,\lambda)\le\mu(p,\theta,\lambda)\le\lambda\quad\forall\,\lambda>0\,.
\ee
If $\theta\in(\theta_\star,1)$, with $q(p,\theta)$ and $\Lambda(p,\theta,\lambda)$ given respectively by~\eqref{qptheta} and~\eqref{Lambda}, then
\be{C2}
\(1-\tfrac p2\,(1-\theta)\)\,\overline{\mu}\(q(p,\theta),\Lambda(p,\theta,\lambda)\)\le\mu(p,\theta,\lambda)\le\lambda\quad\forall\,\lambda>0
\ee
under the condition
\be{Restrictedtheta}
\begin{aligned}
&d=1\,,\quad p\in(2,\infty)\quad\mbox{and}\quad\theta\in(2\,\theta_\star,1)\,,\\
&d=2\,,\quad p\in(2,\infty)\quad\mbox{and}\quad\theta\in(\theta_\star,1)\,,\\
&d\ge3\,,\quad p\in(2,2^*]\quad\mbox{and}\quad\theta\in[\theta_\star,1)\,.
\end{aligned}
\ee
If $\theta\in(\theta_\star,1]$ and $\theta<d/2$, then for any $\lambda>0$, we have
\be{C3}
\frac\lambda\Lambda\,\mu\(\frac{2\,d}{d-2\,\theta},\theta,\Lambda\)\le\mu(p,\theta_\star,\lambda)\le\lambda\quad\mbox{with}\quad\Lambda=\Lambda\(p,\frac{\theta_\star}\theta,\lambda\)=\frac{4\,\theta\,\lambda}{(p-2)\,(d-2\,\theta)}\,.
\ee
\end{corollary}
If $d\ge3$, the lower estimate in~\eqref{C3} is remarkable because $\mu(2^*,1,\Lambda)=\overline{\mu}(2^*,\Lambda)=\min\{\Lambda,d\}$ according, for instance, to Lemma~5 in~\cite{DolEsLa-APDE2014}, so that, using the monotonicity of $\lambda\mapsto\mu(p,\theta_\star,\lambda)$, we obtain
\be{A}
\mu(p,\theta_\star,\lambda)\ge\min\left\{\lambda,(p-2)\,\mathsf A\right\}\quad\mbox{with}\quad\mathsf A:=\frac14\,d\,(d-2)\,.
\ee
{\bf Proof of Corollary~\ref{Cor:mu}.} The upper bounds follow from Lemma~\ref{Lem:mu}, while the lower bounds are consequences of~\eqref{Mon2} with $1=\theta_1\le\theta_2=\theta$ in~\eqref{C1}, $\theta=\theta_1<\theta_2\to1$ in~\eqref{C2}, and $\theta_\star=\theta_1<\theta_2=\theta<1$ in~\eqref{C3}. In~\eqref{C3}, the case $\theta=1$ is covered if $d\ge3$ by taking the limit as $\theta\to1_-$.\hfill\ \qed

As a consequence, we obtain a rigidity range. 
\begin{lemma}\label{Lem:lambda0} Let $d\ge3$. Under Condition~\eqref{Restrictedtheta} and with the above notation, we have
\be{lambda0}
\mu(p,\theta,\lambda)=\lambda\quad\mbox{if}\quad\lambda\le\lambda_0\quad\mbox{and}\quad\mu(p,\theta,\lambda)\ge\lambda_0\quad\mbox{if}\quad\lambda\ge\lambda_0
\ee
for some $\lambda_0\ge d\,(2-p\,(1-\theta))/2\ge(p-2)\,\mathsf A$. Moreover, there is an explicit constant $\kappa_1(p,\theta)>0$ such that
\[
\mu(p,\theta,\lambda)\ge\max\left\{\lambda_0,\kappa_1(p,\theta)\(\tfrac\lambda{(p-2)\,\mathsf A}-1\)^{1-\theta_\star/\theta}\right\}\quad\mbox{if}\quad\lambda\ge\lambda_0\,.
\]
\end{lemma}
\begin{proof} It is known~\cite{1991,DEKL2014,DolEsLa-APDE2014} that $\overline{\mu}(q,\Lambda)=\Lambda$ if $\Lambda\le d$ and $\overline{\mu}(q,\Lambda)\ge d$ if $\Lambda\ge d$. Under Condition~\eqref{Restrictedtheta}, we deduce~\eqref{lambda0} from the lower estimate on $\mu(p,\theta,\lambda)$ in~\eqref{C2} with $\Lambda=\Lambda(p,\theta,\lambda)$ defined by~\eqref{Lambda}. Assume that $d\ge3$ so that $\mathsf A>0$, and let $\mathsf x:=\nrmS u{2^*}^2/\nrmS u2^2$. Using~\eqref{GNS} in the critical case and~\eqref{Holder:p-q} with $q=2^*$, we find that
\[
\mathcal Q_{p,\theta,\lambda}[u]\ge\big((p-2)\,\mathsf A\,\mathsf x+\lambda-(p-2)\,\mathsf A\big)\,\mathsf x^{-\,\theta_\star/\theta}\,.
\]
Assume that $\theta>\theta_\star$. An optimization with respect to $\mathsf x>0$ shows that
\[
\mathcal Q_{p,\theta,\lambda}[u]\ge\kappa_1(p,\theta)\(\frac\lambda{(p-2)\,\mathsf A}-1\)^{1-\theta_\star/\theta}\quad\mbox{where}\quad\kappa_1(p,\theta):=\frac{\theta\,(p-2)\,\mathsf A}{\theta_\star^{\theta_\star/\theta}\,(\theta-\theta_\star)^{1-\theta_\star/\theta}}\,.
\]
The case $\theta=\theta_\star$ is obtained by observing that $\lim_{\theta\to\theta_\star}\kappa_1(p,\theta)=(p-2)\,\mathsf A>0$.\end{proof}

If $d\ge3$, Lemma~\ref{Lem:lambda0} holds with $\lambda_0=(p-2)\,\mathsf A$. The case $d=2$, $\theta=\theta_\star$ and the case $d=1$, $\theta_\star\le\theta\le2\,\theta_\star$ are not covered by Condition~\eqref{Restrictedtheta}. In these cases, lower estimates require a more sophisticated approach that we describe next.

\subsection{Stereographic projection, scalings and interpolation inequalities}\label{Sec:Stereo}

Let us explain how we define and use the stereographic projection. On $\R^d\ni x$, let $r=|x|$ and $\omega=x/|x|$ denote spherical coordinates. On the unit sphere $\S^d\subset\R^{d+1}$, we consider cylindrical coordinates $(\rho\,\omega,z)\in\R^d\times(-1,1)$ with $\rho^2+z^2=1$. The \emph{stereographic projection} $\mathsf S:\S^d\setminus\{\mathsf N\}\to\R^d$, where $\mathsf N\in\S^d$ is the North Pole defined by $z=+1$, is such that $\mathsf S(\rho\,\omega,z)=r\,\omega=x$ where
\[
z=\frac{r^2-1}{r^2+1}\quad\mbox{and}\quad\rho=\frac{2\,r}{1+r^2}\,.
\]
If $v$ is a function on $\R^d$, let us consider its counterpart $u:=\mathsf S^{-1}\,v$ on $\S^d$ obtained using the inverse stereographic projection as
\[
\(u\circ\mathsf S^{-1}\)(x)=\mathsf m(r)^{d-2}\,v(x)\quad\forall\,x\in\R^d\quad\mbox{with}\quad\mathsf m(r)=\sqrt{(1+r^2)/2}\,.
\]
For any $q\ge1$, with $\delta(q):=2\,d-q\,(d-2)$, we obtain
\begin{align*}
&\isd{|u|^q}=\big|\S^d\big|^{-1}\ird{\frac{|v|^q}{\mathsf m(r)^{\delta(q)}}}\,,\\
&\isd{|\nabla u|^2}+\frac 14\,d\,(d-2)\isd{|u|^2}=\big|\S^d\big|^{-1}\ird{|\nabla v|^2}\,.
\end{align*}
Notice that $\delta(2)=4$, and $\delta(2^*)=0$ if $d\ge3$. We refer to Theorem~2.1 of~\cite{Dolbeault_2020b} for a statement concerning the inequality obtained from~\eqref{Ineq:GNS1} and a special case of~\eqref{Ineq:GNS2} by the stereographic projection.

Using the stereographic projection, we learn that
\[
\mu(p,\theta,\lambda)=|\S^d|^{-\frac{p-2}{p\,\theta}}\,\inf\,\frac{(p-2)\,\nrmR{\nabla v}2^2+4\,(\lambda-\mathsf A)\ird{\frac{|v|^2}{(1+|x|^2)^2}}}{\(\ird{\big(\frac2{1+|x|^2}\big)^{\delta/2}\,|v|^p}\)^{2/(p\,\theta)}\(\ird{\frac{4\,|v|^2}{(1+|x|^2)^2}}\)^{1-1/\theta}}
\]
where the infimum is taken on $\left\{v\in\mathrm L^2\big(\R^d,(1+|x|^2)^{-2}\,dx\big)\setminus\{0\}\,:\,
\nabla v\in\mathrm L^2(\R^d,dx)\right\}$ and $\delta=\delta(p)$. Either $v$ is constant, in which case $\mathcal Q_{p,\theta,\lambda}[u]=\lambda$, or $\nrmR{\nabla v}2>0$, which we assume from now on. With the rescaling
\[
v(x)=w\(\frac x\tau\)\quad\forall\,x\in\R^d\,,\quad\tau=\frac1{\sqrt{\lambda-\mathsf A}}
\]
for any $\lambda>\mathsf A$, we find that
\begin{multline*}
\mathcal Q_{p,\theta,\lambda}[u]:=|\S^d|^{-\frac{2\,\theta_\star}{d\,\theta}}\,\frac{(p-2)\,\nrmR{\nabla v}2^2+4\,(\lambda-\mathsf A)\ird{\frac{|v|^2}{(1+|x|^2)^2}}}{\(\ird{\big(\frac2{1+|x|^2}\big)^{\delta/2}\,|v|^p}\)^{2/(p\,\theta)}\(\ird{\frac{4\,|v|^2}{(1+|x|^2)^2}}\)^{1-1/\theta}}\\
=\(\frac\tau{|\S^d|^{1/d}}\)^\frac{2\,\theta_\star}\theta\,\frac{\tau^{-2}\,(p-2)\,\nrmR{\nabla w}2^2+4\,(\lambda-\mathsf A)\ird{\frac{|w|^2}{(1+\tau^2\,|x|^2)^2}}}{\(\ird{\big(\frac2{1+\tau^2\,|x|^2}\big)^{\delta/2}\,|w|^p}\)^{2/(p\,\theta)}\(\ird{\frac{4\,|w|^2}{(1+\tau^2\,|x|^2)^2}}\)^{1-1/\theta}}\\
\hspace*{2cm}=\frac{(\lambda-\mathsf A)^{1-\theta_\star/\theta}}{|\S^d|^{2\,\theta_\star/(d\,\theta)}}\,\frac{(p-2)\,\nrmR{\nabla w}2^2+4\ird{\frac{|w|^2}{(1+\tau^2\,|x|^2)^2}}}{\(\ird{\big(\frac2{1+\tau^2\,|x|^2}\big)^{\delta/2}\,|w|^p}\)^{2/(p\,\theta)}\(\ird{\frac{4\,|w|^2}{(1+\tau^2\,|x|^2)^2}}\)^{1-1/\theta}}\,.
\end{multline*}
Assuming that $X$ and $Y$ are positive numbers and using the identity
\be{ctheta}
\(\frac{X+Y}{Y^{1-1/\theta}}\)^\theta\ge\mathsf c_\theta\,X^{\theta_\star}\,Y^{1-\theta_\star}\quad\mbox{where}\quad\mathsf c_\theta:=\frac{\theta^\theta}{\theta_\star^{\,\theta_\star}\,\,(\theta-\theta_\star)^{\theta-\theta_\star}}
\ee
with $X=(p-2)\,\nrmR{\nabla w}2^2$ and $Y=4\ird{\frac{|w|^2}{(1+\tau^2\,|x|^2)^2}}$, we obtain
\[
\mathcal Q_{p,\theta,\lambda}[u]\ge\frac{(p-2)^{\theta_\star/\theta}\,(\lambda-\mathsf A)^{1-\theta_\star/\theta}}{|\S^d|^{2\,\theta_\star/(d\,\theta)}}\(c_\theta\,\frac{\(\ird{|\nabla w|^2}\)^{\theta_\star}\(\ird{\frac{|w|^2}{(1+\tau^2\,|x|^2)^2}}\)^{1-\theta_\star}}{\(\ird{\big(\frac{|w|^p}{(1+\tau^2\,|x|^2)^{\delta/2}}}\)^{2/p}}\)^{1/\theta}\,.
\]
By yet another scaling, we can get rid of $(1+\tau^2\,|x|^2)$ and find that, after minimization,
\be{h}
\mathcal Q_{p,\theta,\lambda}[u]\ge\kappa_2(p,\theta)\,(\lambda-\mathsf A)^{1-\theta_\star/\theta}\quad\mbox{with}\quad\kappa_2(p,\theta):=(p-2)^{\theta_\star/\theta}\,{c_\theta^{1/\theta}}\(\tfrac{\mathcal C_{\mathrm{GNS}}(p)}{|\S^d|^{\theta_\star/d}}\)^{2/\theta}\,,
\ee
where $\mathcal C_{\mathrm{GNS}}(p)$ is the optimal constant in~\eqref{Ineq:Rd}.
\begin{proposition}\label{Prop:lambda0} Let $d\ge1$. If $p>2$, $p\le2^*$ if $d\ge3$, and $\theta\ge\theta_\star$, with $\kappa_2(p,\theta)$ defined by~\eqref{h}, we have
\[
\mu(p,\theta,\lambda)\ge\min\left\{\lambda,\kappa_2(p,\theta)\,(\lambda-\mathsf A)^{1-\theta_\star/\theta}\right\}\quad\forall\,\lambda>\mathsf A\,.
\]
Moreover, if $d\le2$, there is some $\lambda_0>0$ such that $\mu(p,\theta,\lambda)=\lambda$ if $\lambda\le\lambda_0$.
\end{proposition}
\begin{proof} Either the optimal function is a constant, in which case \hbox{$\mu(p,\theta,\lambda)=\lambda$}, or~\eqref{h} applies. If $d\le2$, we have $\mathsf A\le0$ which gives an explicit $\lambda_0>0$, obtained by solving $\lambda_0=\kappa_2(p,\theta)\,(\lambda_0-\mathsf A)^{1-\theta_\star/\theta}$. As a special case, we have $\lambda_0=(\kappa_2(p,\theta))^{\theta/\theta_\star}$ if $d=2$.\end{proof}

\subsection{Existence of an optimal function}

\begin{lemma}\label{Lem:Existence} Let $d\ge1$, $p\in(2,+\infty)$ if $d=1$, $2$, $p\in(2,2^*]$ if $d\ge3$, and $\theta>\theta_\star$. For any $\lambda>0$, there is an optimal function for~\eqref{Ineq:GNS2} if either $d\le2$, or $d\ge3$ and~$p<2^*$.\end{lemma}
\begin{proof} If $d\ge3$, H\"older's inequality~\eqref{Holder:p-q} applies with $\theta=\theta_\star$ and $q(p,\theta_\star)=2^*$. Combined with~\eqref{GNS} written with $p=2^*$, we obtain
\[
\nrmS up\le\nrmS u{2^*}^{\theta_\star}\,\nrmS u2^{1-\theta_\star}\le\(\frac4{d\,(d-2)}\nrmS{\nabla u}2^2+\nrmS u2^2\)^{\theta_\star/2}\,\nrmS u2^{1-\theta_\star},
\]
and, as a consequence, using~\eqref{mu-Q},
\[
\mathcal Q_{p,\theta,\lambda}[u]\ge\frac{(p-2)\,\nu+\lambda}{\(\frac4{d\,(d-2)}\,\nu+1\)^{\theta_\star/\theta}}=:\mathsf q_{p,\theta,\lambda}[\nu]\quad\mbox{with}\quad\nu=\nu[u]:=\frac{\nrmS{\nabla u}2^2}{\nrmS u2^2}\,.
\]
The existence of an optimal function follows by a standard minimization argument. If $(u_n)_{n\in\N}$ is a minimizing sequence for $\mathcal Q_{p,\theta,\lambda}$, we can assume without loss of generality that $\nrmS{u_n}2=1$ by homogeneity of $\mathcal Q_{p,\theta,\lambda}$. Using $\theta>\theta_\star$, it follows from the above estimate that $\nu[u_n]$ is bounded, thus proving that $(u_n)_{n\in\N}$ weakly converges in $\mathrm H^1(\S^d)$ up to the extraction of a subsequence to some $u\in\mathrm H^1(\S^d)$. From the compactness of the embedding $\mathrm H^1(\S^d)\hookrightarrow\mathrm L^p(\S^d)$, we conclude that $u$ is a minimizer of $\mathcal Q_{p,\theta,\lambda}$. In dimension $d=1$ and $d=2$, a similar method applies except that one has to interpolate with an $\mathrm L^q(\S^d)$ norm for some $q\in(p,2^*)=(p,+\infty)$ instead of picking $q=2^*$.\end{proof}

Notice that an optimization of $\nu\mapsto\mathsf q_{p,\theta,\lambda}[\nu]$ with respect to $\nu\in\R^+$ gives an explicit lower estimate of $\mathcal Q_{p,\theta,\lambda}[u]$ which is similar to the optimization on $\mathsf x$ in the proof of Lemma~\ref{Lem:lambda0}.

\subsection{Symmetry breaking}\label{Sec:SymmetryBreaking}

We obtain an estimate of the symmetry breaking range by Taylor expanding around a constant perturbed by a spherical harmonic function.
\begin{lemma}\label{Lem:SymmetryBreaking} Let $d\ge1$, $p\in(2,2^*)$ if $d=1$, $2$, $p\in(2,2^*]$ if $d\ge3$ and $\theta>\theta_\star$. There is \emph{symmetry breaking} in~\eqref{Ineq:GNS2}, that is, $\mu(p,\theta,\lambda)<\lambda$, if $\lambda>d\,\theta$.\end{lemma}
\begin{proof} Let us consider $u_\varepsilon=1+\varepsilon\,\varphi$ with $\varphi\in\mathrm H^1(\S^d,d\sigma)$ such that $\isd\varphi=0$ as a test function:
\begin{align*}
&\nrm{\nabla u_\varepsilon}2^2=\varepsilon^2\,\nrmS{\nabla\varphi}2^2\,,\\
&\nrm{u_\varepsilon}2^2=1+\varepsilon^2\,\nrmS\varphi2^2\,,\\
&\nrm{u_\varepsilon}p^2=1+\varepsilon^2\,(p-1)\,\nrmS\varphi2^2+o(\varepsilon^2)\,.
\end{align*}
As a consequence, we obtain that
\[
\lim_{\varepsilon\to0}\frac1{\varepsilon^2}\,\frac{\mathcal Q_{p,\theta,\lambda}[u_\varepsilon]-\lambda}{p-2}=\nrmS{\nabla\varphi}2^2-\frac\lambda\theta\,\nrmS\varphi2^2=\(d-\frac\lambda\theta\)\nrmS\varphi2^2<0
\]
if we take a non-trivial spherical harmonic function $\varphi$ such that $-\,\Delta\varphi=d\,\varphi$ and $\lambda>d\,\theta$.\end{proof}

\subsection{Simple symmetry results}

Using a simple monotonicity observation, we obtain the following result.
\begin{lemma}\label{Lem:Symmetry} Let $d\ge1$, $p\in(2,2^\#)$. If $\theta\ge\theta^\#$, there is \emph{symmetry} in~\eqref{Ineq:GNS2}, that is, $\mu(p,\theta,\lambda)=\lambda$, if $\lambda\le d\,\theta$.\end{lemma}
\begin{proof} We start by the following observation. Let $X>0$ be a given real number and consider the function $\mathsf f$ on $\R^+$ defined by
\[
\mathsf f(\theta):=\theta\(X^\frac1\theta-1\)\,.
\]
It is elementary to compute $-\,\mathsf f'(\theta)=X^{\frac1\theta}\,\log(X^{\frac1\theta})+1-X^{\frac1\theta}\ge0$.

Let $X:=\nrmS up^2/\nrmS u2^2$. If $\lambda=d\,\theta$, we consider two cases:\\
(1) Either $\theta\ge1$ and we read from~\eqref{GNS} that
\begin{align*}
\nrmS{\nabla u}2^2&\ge\;\frac d{p-2}\(\nrmS up^2-\nrmS u2^2\)=\frac d{p-2}\,\nrmS u2^2\,\mathsf f(1)\\
&\ge\;\frac d{p-2}\,\nrmS u2^2\,\mathsf f(\theta)=\frac{d\,\theta}{p-2}\(\nrmS up^\frac2\theta\,\nrmS u2^{-2\,\frac{1-\theta}\theta}-\nrmS u2^2\)\,.
\end{align*}
(2) Or $\theta\in(\theta^\#,1)$ and we know from~\eqref{Ineq:GNSimproved} that
\begin{align*}
\nrmS{\nabla u}2^2&\ge\;\frac d{p-2}\,\nrmS u2^2\,\mathsf f(\theta^\#)\\
&\ge\;\frac d{p-2}\,\nrmS u2^2\,\mathsf f(\theta)=\frac{d\,\theta}{p-2}\(\nrmS up^\frac2\theta\,\nrmS u2^{-2\,\frac{1-\theta}\theta}-\nrmS u2^2\)\,.
\end{align*}
If $\lambda<d\,\theta$, the result follows from
\begin{align*}
\nrmS{\nabla u}2^2\ge&\;\frac{d\,\theta}{p-2}\(\nrmS up^\frac2\theta\,\nrmS u2^{-2\,\frac{1-\theta}\theta}-\nrmS u2^2\)\\
&\;=\frac\lambda{p-2}\(\nrmS up^\frac2\theta\,\nrmS u2^{-2\,\frac{1-\theta}\theta}-\nrmS u2^2\)+\frac{d\,\theta-\lambda}{p-2}\(\nrmS up^\frac2\theta\,\nrmS u2^{-2\,\frac{1-\theta}\theta}-\nrmS u2^2\)\\
&\ge\;\frac\lambda{p-2}\(\nrmS up^\frac2\theta\,\nrmS u2^{-2\,\frac{1-\theta}\theta}-\nrmS u2^2\)\,,
\end{align*}
as a consequence of $\nrmS up^{2/\theta}\,\nrmS u2^{-2\,(1-\theta)/\theta}\ge\nrmS u2^2$, \emph{i.e.}, of H\"older's inequality $\nrmS u2\le\nrmS up$. The equality case means that $u$ is a constant.\end{proof}

\subsection{Consequences for the main results}\label{Sec:Main1}

Let us collect the consequences of the above estimates for Theorems~\ref{Thm1} and~\ref{Thm2}. The proofs (part 2/2) will be completed in Section~\ref{Sec:Main2}.

\medskip\noindent{\bf Proof of Theorem~\ref{Thm1}, part 1/2.} Inequality~\eqref{Ineq:GNS2} holds as a consequence of Corollary~\ref{Cor:mu} if $\theta>(p-2)/p$ or equivalently Condition~\eqref{Restrictedtheta}, and of Proposition~\ref{Prop:lambda0} if either $d=2$ and $\theta=\theta_\star$, or $d=1$ and $\theta\in[\theta_\star,2\,\theta_\star]$. The existence of optimal functions in the equality cases of~\eqref{Ineq:GNS2} is established in Lemma~\ref{Lem:Existence} while the monotonicity and the concavity of $\lambda\mapsto\mu(p,\theta,\lambda)$ is a straightforward consequence of the expression of $\mathcal Q_{p,\theta,\lambda}$ given by~\eqref{mu-Q}. The upper bound $\mu\le\lambda$ and the symmetry breaking property meaning $\mu<\lambda$ are the result of Lemma~\ref{Lem:mu} and Lemma~\ref{Lem:SymmetryBreaking} respectively. The symmetry result valid for $p<2^\#$ and $\theta\ge\theta^\#$ corresponds to Lemma~\eqref{Lem:Symmetry}. As a consequence, the proof is complete except for the proof of~\eqref{asymptotics}.\hfill\ $\square$

\medskip\noindent{\bf Proof of Theorem~\ref{Thm2}, part 1/2.} We deduce that $\lambda_{\mathrm T}$ is positive from Lemma~\ref{Lem:lambda0} if $d\ge3$ and from Proposition~\ref{Prop:lambda0} if $d=1$ or $d=2$. The fact that $\lambda_{\mathrm T}(p,\theta)\le d\,\theta$ is a consequence of Theorem~\ref{Thm1} (or Lemma~\ref{Lem:SymmetryBreaking}) but the strict inequality for $\theta-\theta_\star\ge0$ small enough and $p$ in the  appropriate range is still to be proved.\hfill\ $\square$

\section{Asymptotics, reparametrization and further estimates}\label{Sec:More}

\subsection{Asymptotics of \texorpdfstring{$\mu(p,\theta,\lambda)$ as $\lambda\to+\infty$}{}}

Heuristically, large values of $\lambda$ correspond to highly concentrated optimal functions, which behave like the functions obtained by a \emph{stereographic projection} on the tangent Euclidean space at the point of concentration as in Section~\ref{Sec:Stereo}. Here we shall use the optimal constant $\mathcal K_{p,d,\theta}$ in the Gagliardo-Ni\-ren\-berg-Sobolev inequality
\be{Ineq:Rd-inhom}
\(\nrmR{\nabla v}2^2+\nrmR v2^2\)\nrmR v2^{\frac2\theta-2}\ge\mathcal K_{p,d,\theta}\,\nrmR vp^\frac2\theta\quad\forall\,v\in\mathrm H^1(\R^d,dx)\,.
\ee
\begin{lemma}\label{Lem:K} Let $d\ge1$, $p>2$, $p\le2^*$ if $d\ge3$, and $\theta\ge\theta_\star$. With $c_\theta$ as in~\eqref{ctheta}, we have
\[
\mathcal K_{p,d,\theta}=\big(c_\theta\,\mathcal C_{\mathrm{GNS}(p)}^2\big)^{1/\theta}\,.
\]\end{lemma}
\noindent Inequality~\eqref{Ineq:Rd-inhom} is equivalent to~\eqref{Ineq:Rd}, as will be made clear in the proof.
\begin{proof} Assume that $\theta>\theta_\star$ and let us consider $u_\varepsilon(x):=\varepsilon^{-d/p}\,u(x/\varepsilon)$ for any $x\in\R^d$, 
so that
\[
\ird{|u_\varepsilon|^p}=\ird{|u|^p}\,,\quad\ird{|u_\varepsilon|^2}=\varepsilon^{2\,\theta_\star}\ird{|u|^2}\quad\mbox{and}\quad\ird{|\nabla u_\varepsilon|^2}=\varepsilon^{-2\,(1-\theta_\star)}\ird{|\nabla u|^2}\,.
\]
For an arbitrary $\lambda>0$, we compute
\begin{multline*}
\(\nrmR{\nabla u_\varepsilon}2^2+\frac\lambda{p-2}\,\nrmR{u_\varepsilon}2^2\)^\theta\Bigg(\nrmR{u_\varepsilon}2^2\Bigg)^{1-\theta}\\
=\(\varepsilon^{-2\,(1-\theta_\star)}\,\nrmR{\nabla u}2^2+\frac\lambda{p-2}\,\varepsilon^{2\,\theta_\star}\,\nrmR u2^2\)^\theta\Bigg(\varepsilon^{2\,\theta_\star}\,\nrmR u2^2\Bigg)^{1-\theta}=\(\mathsf A\,\varepsilon^{-\mathsf a}+\frac{\lambda\,\mathsf B}{p-2}\,\varepsilon^{\mathsf b}\)^\theta
\end{multline*}
where $\mathsf A:=\nrmR{\nabla u}2^2\,\nrmR u2^{2\,(1-\theta)/\theta}$, $\mathsf B:=\nrmR u2^{2/\theta}$, $\mathsf a:=2-\mathsf b$ and $\mathsf b:=2\,\theta_\star/\,\theta$. With $c_\theta$ defined as in~\eqref{ctheta}, an optimization on $\varepsilon>0$ shows that
\be{K-CGNS}
\(\mathsf A\,\varepsilon^{-\mathsf a}+\frac{\lambda\,\mathsf B}{p-2}\,\varepsilon^{\mathsf b}\)^\theta=c_\theta\(\frac\lambda{p-2}\)^{\theta-\theta_\star}\,\nrmR{\nabla u}2^{2\,\theta_\star}\,\nrmR u2^{2\,(1-\theta_\star)}\ge c_\theta\(\frac\lambda{p-2}\)^{\theta-\theta_\star}\,\mathcal C_{\mathrm{GNS}}(p)^2\,\nrmR up^2
\ee
using~\eqref{Ineq:Rd}. Taking $\lambda=p-2$, we find that $\mathcal K_{p,d,\theta}\ge c_\theta^{1/\theta}\,\mathcal C_{\mathrm{GNS}(p)}^{2/\theta}$. There is in fact equality because we can choose for $u$ an optimal function for~\eqref{Ineq:Rd}. The case $\theta=\theta_\star$ is obtained as a limit case using $\displaystyle\lim_{\theta\to\theta_\star}c_\theta=1$.\end{proof}

The above computation is also interesting if $\lambda\neq p-2$. On $\S^d$, in the limit as $\lambda\to+\infty$, so that \hbox{$\varepsilon\to0$}, heuristics are simple as soon as we know that optimal functions concentrate: the behaviour of $\lambda\mapsto\mu(p,\theta,\lambda)$ as $\lambda\to+\infty$ is determined by the interpolation problem on the tangent hyperplane to the sphere at the point of concentration, \emph{i.e.}, by~\eqref{Ineq:Rd-inhom}, up to the normalization factor in the definition of $d\sigma$. This was proved for $\theta=1$ in~\cite{DolEsLa-APDE2014}. By taking for $u$ an approximation of the optimal function for~\eqref{Ineq:Rd} with compact support in a small ball, which can be done because of the scaling invariance of~\eqref{Ineq:Rd}, we can guess from~\eqref{K-CGNS} that $\mu(p,\theta,\lambda)=O(\lambda^{\theta-\theta_\star})$ as $\lambda\to+\infty$. The purpose of the following result is to extend the results for~\eqref{Ineq:GNS1} of Proposition~9 and~10 of~\cite{DolEsLa-APDE2014} to Inequality~\eqref{Ineq:GNS2} with $\theta\neq1$.
\begin{lemma}\label{Lem:Asymptotics} Under the assumptions of Theorem~\ref{Thm1} on $d$, $p$ and $\theta\ge\theta_\star$, we have
\[
\mu(p,\theta,\lambda)=\kappa(p,\theta)\,\lambda^{1-\theta_\star/\theta}\,\big(1+o(1)\big)\quad\mbox{as}\quad\lambda\to+\infty\quad\mbox{where}\quad\kappa(p,\theta)=\(\tfrac{p-2}{|\S^d|^{2/d}}\)^{\theta_\star/\theta}\mathcal K_{p,d,\theta}\,.
\]
\end{lemma}
\begin{proof} We do a \emph{blow-up analysis} based on a change of variables depending on $\lambda\to+\infty$. As in Section~\ref{Sec:Stereo}, let $\mathsf m(r)=\sqrt{(1+r^2)/2}$, $\delta(p)=2\,d-p\,(d-2)$ and $u_\lambda:=\mathsf S^{-1}\,v_\lambda$ where
\[
v_\lambda(x):=v\(\frac x\tau\)\quad\mbox{with}\quad\tau=\frac12\,\sqrt{\frac{p-2}{\lambda-\overline{\lambda}}}\quad\mbox{and}\quad\overline{\lambda}:=(p-2)\,\mathsf A=\frac14\,d\,(d-2)\,(p-2)\,.
\]
We take $v$ as an optimal function for~\eqref{Ineq:Rd-inhom}. With $r=|x|$ and $\mathcal Q_{p,\theta,\lambda}$ defined in~\eqref{mu-Q}, we compute
\begin{align*}
\mathcal Q_{p,\theta,\lambda}[u_\lambda]&=\frac{(p-2)\(\nrmS{\nabla u_\lambda}2^2+\frac 14\,d\,(d-2)\,\nrmS{u_\lambda}2^2\)+(\lambda-\overline{\lambda})\,\nrmS{u_\lambda}2^2}{\nrmS{u_\lambda}p^{\frac2\theta}\,\nrmS{u_\lambda}2^{2-\frac2\theta}}\\
&=\frac{(p-2)\,\nrmR{\nabla v_\lambda}2^2+(\lambda-\overline{\lambda})\,\nrmR{\mathsf m(r)^{-2}\,v_\lambda}2^2}{|\S^d|^{\frac1\theta\(1-\frac2p\)}\,\nrmR{\mathsf m(r)^{-\frac{\delta(p)}p}\,v_\lambda}p^{\frac2\theta}\,\nrmR{\mathsf m(r)^{-2}\,v_\lambda}2^{2-\frac2\theta}}\\
&=\frac{(\lambda-\overline{\lambda})^{1-\frac{\theta_\star}\theta}\,(p-2)^{\frac{\theta_\star}\theta}\,\(\nrmR{\nabla v}2^2+\frac14\,\nrmR{\mathsf m(\tau\,r)^{-2}\,v}2^2\)}{|\S^d|^{\frac1\theta\(1-\frac2p\)}\,2^{-\frac{\delta(p)}{p\,\theta}-2\left(1-\frac1\theta\right)}\,\nrmR{\mathsf m(\tau\,r)^{-\frac{\delta(p)}p}\,v}p^{\frac2\theta}\,\nrmR{\mathsf m(\tau\,r)^{-2}\,v}2^{2-\frac2\theta}}\,.
\end{align*}
By taking the limit as $\lambda\to+\infty$ and using Lebesgue's theorem of dominated convergence with $\mathsf m(\tau\,r)\to1/\sqrt2$ a.e., and using the fact that the quotient becomes proportional to the quotient that defines the constant in the Gagliardo-Nirenberg-Sobolev inequality up to a numerical constant, we conclude that
\[
\mathcal Q_{p,\theta,\lambda}[u_\lambda]=(\lambda-\overline{\lambda})^{1-\theta_\star/\theta}\kappa(p,\theta)\big(1+o(1)\big)
\]
as $\lambda\to+\infty$, for a constant $\kappa(p,\theta)$ that can be computed explicitly, up to tedious but elementary considerations. Consistently with Proposition~\ref{Prop:lambda0} (we notice that $\kappa(p,\theta)=\kappa_2(p,\theta)$ in~\eqref{h}), this proves that 
\[
\lim_{\lambda\to+\infty}\lambda^{\theta_\star/\theta-1}\,\mu(p,\theta,\lambda)\le\kappa(p,\theta).
\]

We have now to prove the reverse inequality. As in Proposition~10 of~\cite{DolEsLa-APDE2014}, we argue by contradiction. Let $\seq\lambda n$ and $\seq un$ be such that $u_n\in\mathrm H^1(\S^d,d\sigma)$ with $\nrmS{u_n}p=1$ for any $n\in\N$, $\lim_{n\to+\infty}\lambda_n=+\infty$,
\be{splitting}
\mathcal Q_{p,\theta,\lambda_n}[u_n]=\mu(p,\theta,\lambda_n)\quad\mbox{and}\quad\lim_{n\to+\infty}\lambda_n^{\theta_\star/\theta-1}\,\mu(p,\theta,\lambda_n)\le\kappa(p,\theta)-\eta
\ee
for some $\eta>0$. We learn from the expression of $\mathcal Q_{p,\theta,\lambda}$ that $\lambda_n\,\nrmS{u_n}2^{2/\theta}\le\kappa(p,\theta)\,\lambda_n^{1-\theta_\star/\,\theta}$ as $n\to+\infty$ so that $\lim_{n\to+\infty}\nrmS{u_n}2=0$. Concentration occurs: there exists a sequence $\seq yi$ of points in $\S^d$ with $i\in I\subset\N$, sequences of positive numbers $\seq\zeta i$ and $(r_{i,n})_{i,n\in I\times\N}$ and functions $u_{i,n}\in\mathrm H^1(\S^d,d\sigma)$ with
\[
u_{i,n} =u_n\quad\mbox{on}\quad\S^d\cap B(y_i,r_{i,n})\quad\mbox{and}\quad\mbox{supp }u_{i,n}\subset \S^d\cap B(y_i,2\,r_{i,n})
\]
such that, as $n\to+\infty$,
\[
\sum_{i\in I}\nrmS{\nabla u_{i,n}}2^2\sim\nrmS{\nabla u_n}2^2\quad\mbox{and}\quad\sum_{i\in I}\nrmS{u_{i,n}}2^2\sim\nrmS{u_n}2^2
\]
with, for all $i\in I$,
\[
\lim_{n\to+\infty}r_{i,n}=0\,,\quad\sum_{i\in\N}\,\zeta_i=1\ \quad\mbox{and}\quad\lim_{n\to+\infty}\int_{\S^d\cap B(y_i,r_{i,n})}|u_{i,n}|^p\,d\sigma=\zeta_i\,.
\]
With a \emph{blow up} argument applied to $(u_{i,n})_{n\in\N}$, we prove for all $i\in I$ that
\[
\lim_{n\to+\infty}\lambda_n^{\frac{\theta_\star}\theta-1}\((p-2)\,\nrmS{\nabla u_{i,n}}2^2+\lambda_n\,\nrmS{u_{i,n}}2^2\)\nrmS{u_{i,n}}2^{-2-\frac2\theta}\ge\kappa(p,\theta)\,\zeta_i^{\frac2{p\,\theta}}\,.
\]
We may notice that
\[
\frac{(p-2)\,\nrmS{\nabla u_{i,n}}2^2+\lambda_n\,\nrmS{u_{i,n}}2^2}{\nrmS{u_{i,n}}2^{2-\frac2\theta}}\ge\frac{(p-2)\,\nrmS{\nabla u_{i,n}}2^2+\lambda_n\,\nrmS{u_{i,n}}2^2}{\nrmS{u_n}2^{2-\frac2\theta}}\,.
\]
Let us choose an integer $N$ such that
\[
\sum_{i=1}^N \zeta_i^{\frac2{p\,\theta}}\ge\(\sum_{i=1}^N \zeta_i\)^{\frac2{p\,\theta}}>1-\frac\eta{2\,\kappa(p,\theta)}\,.
\]
For $n$ large enough, by writing
\[
\lambda_n^{\frac{\theta_\star}\theta-1}\,\sum_{i=1}^N\frac{(p-2)\,\nrmS{\nabla u_{i,n}}2^2+\lambda_n\,\nrmS{u_{i,n}}2^2}{\nrmS{u_{i,n}}2^{2-\frac2\theta}}\ge\kappa(p,\theta)\sum_{i=1}^N\zeta_i^{\frac2{p\,\theta}}\,\big(1+o(1)\big)\ge\kappa(p,\theta)-\frac\eta2\,,
\]
we obtain a contradiction with~\eqref{splitting}. This concludes the proof.
\end{proof}

We learn from Lemma~\ref{Lem:K} and Lemma~\ref{Lem:Asymptotics} that
\be{Bound:thetastar}
\mu(p,\theta_\star,\lambda)\le\kappa(p,\theta_\star)=\tfrac{p-2}{|\S^d|^{2/d}}\,\mathcal C_{\mathrm{GNS}}(p)^{2/\theta_\star}\quad\forall\,\lambda>0\,.
\ee
This estimate is very useful. Taking into account various other results, here is one of the consequences.
\begin{corollary}\label{Cor:LambdaT} With the notation of Theorem~\ref{Thm2}, we have $\lambda_{\mathrm T}(p,\theta)\ge\lambda_{\mathrm T}(p,\theta_\star)$ for any $\theta\ge\theta_\star$ and
\begin{align*}
&(p-2)\,\mathsf A\le\lambda_{\mathrm T}(p,\theta_\star)\le\kappa(p,\theta_\star)\quad\mbox{if}\quad d\ge3\,,\\
&\lambda_{\mathrm T}(p,\theta_\star)=\kappa(p,\theta_\star)\quad\mbox{if}\quad d=1\quad\mbox{or}\quad d=2\,.
\end{align*}
\end{corollary} 
\begin{proof} By Lemma~\ref{Lem:mu} and Lemma~\ref{Lem:theta}, using the monotonicity of $\theta\mapsto\mu(p,\theta,\lambda)$, we know that $\mu(p,\theta_\star,\lambda)\le\mu(p,\theta,\lambda)\le\lambda$ if $\theta\ge\theta_\star$. This proves $\lambda_{\mathrm T}(p,\theta)\ge\lambda_{\mathrm T}(p,\theta_\star)$. If $d\ge3$, the lower estimate on $\lambda_{\mathrm T}(p,\theta_\star)$ follows from~\eqref{A}, while the upper estimate $\lambda_{\mathrm T}(p,\theta_\star)\le\kappa(p,\theta_\star)$ is a consequence of~\eqref{Bound:thetastar}. If $d=1$ or $d=2$, $\lambda_{\mathrm T}(p,\theta_\star)=\kappa(p,\theta_\star)$ is deduced from Proposition~\ref{Prop:lambda0}.\end{proof}

\begin{remark} In view of Corollary~\ref{Cor:LambdaT}, one can conjecture that $\lambda_{\mathrm T}(p,\theta_\star)=\kappa(p,\theta_\star)$ if $d\ge3$, as it is the case if $d=1$ or $d=2$. One the other hand, one may wonder if Lemma~\ref{Lem:Asymptotics} is compatible with~\eqref{C3}, because for a given $p$ and a given $\theta$, $\lambda/\Lambda=(p-2)\,(d-2\,\theta)/4\,\theta$ is independent of $\Lambda>0$. We know by~\eqref{Bound:thetastar} that $\mu(p,\theta_\star,\lambda)$ is bounded uniformly with respect to $\lambda>0$, while $\mu\(q,\theta,\Lambda\)\sim\frac{4\,\theta\,\kappa(q,\theta)}{(p-2)\,(d-2\,\theta)}\,\lambda^{1-\theta_\star(q)/\theta}\to+\infty$ as $\lambda\to+\infty$ if $\theta_\star(q)=d\,(q-2)/(2\,q)>\theta$. There is no paradox because we obtain $\theta_\star(q)=\theta$ if we take $q=2\,d/(d-2\,\theta)$.\end{remark}

\subsection{Reparametrization}\label{Sec:Reparametrization}

As in~\cite{Caffarelli1984,delPino20102045,springerlink:10.1007/s00526-011-0394-y,FreefemDolbeaultEsteban,0951-7715-27-3-435,DoEsFiTe2014},
the set of optimal functions for~\eqref{mu-Q} for $\theta\neq1$ can be parametrized by the critical points of the problem corresponding to $\theta=1$. The precise statement goes as follows.
\begin{proposition}\label{Prop0:Reparametrization} Let $d\ge1$, $p\in(2,2^*)$, $\theta\ge\theta_\star$ and $\Lambda>0$, given. If $u$ is an optimal function for~\eqref{mu-Q}, then, up to a multiplicative constant, the function $u$ solves
\be{EL0}
-\,\Delta u+\frac\lambda{p-2}\,u=u^{p-1}
\ee
with
\be{lambda-mu}
\lambda=\frac1\theta\(\Lambda+(1-\theta)\,(p-2)\,\tfrac{\nrmS{\nabla u}2^2}{\nrmS u2^2}\)\quad\mbox{and}\quad\nrmS up^p=\frac1\theta\(\nrmS{\nabla u}2^2+\frac\Lambda{p-2}\,\nrmS u2^2\)\,.
\ee
\end{proposition}
\begin{proof} Any optimal function $u$ for~\eqref{mu-Q} is such that
\[
\(\nrmS{\nabla u}2^2+\frac\Lambda{p-2}\,\nrmS u2^2\)^\theta\Bigg(\nrmS u2^2\Bigg)^{1-\theta}=\(\frac{\mu(p,\theta,\Lambda)}{p-2}\)^\theta\,\nrmS up^2
\]
solves the Euler-Lagrange equation
\be{EL}
\frac{\theta\,\big(-\,\Delta u+\frac\Lambda{p-2}\,u\big)}{\nrmS{\nabla u}2^2+\frac\Lambda{p-2}\,\nrmS u2^2}+\frac{(1-\theta)\,u}{\nrmS u2^2}=\frac{u^{p-1}}{\nrmS up^p}\,,
\ee
that is, Equation~\eqref{EL0} with $\lambda$ and the multiplicative constant given by~\eqref{lambda-mu}.
\end{proof}

An interesting consequence is the fact that the curve $\Lambda\mapsto\mu(p,\theta,\Lambda)$ can be parametrized by $\lambda$ using
\[
\Lambda=\theta\,\lambda-(1-\theta)\,(p-2)\,\frac{\nrmS{\nabla u}2^2}{\nrmS u2^2}\quad\mbox{and}\quad\mu(p,\theta,\Lambda)=\theta\(\lambda+(p-2)\,\frac{\nrmS{\nabla u}2^2}{\nrmS u2^2}\)\(\frac{\nrmS u2}{\nrmS up}\)^{2/\theta}\,,
\]
where $u$ is as in the above proof, because
\[
\((p-2)\,\frac{\nrmS{\nabla u}2^2}{\nrmS u2^2}+\Lambda\)^\theta\,\frac{\nrmS u2^2}{\nrmS up^2}=\big(\mu(p,\theta,\Lambda)\big)^\theta\;\mbox{and}\;(p-2)\,\frac{\nrmS{\nabla u}2^2}{\nrmS u2^2}+\Lambda=\theta\(\lambda+(p-2)\,\frac{\nrmS{\nabla u}2^2}{\nrmS u2^2}\).
\]
However, it is not \emph{a priori} granted that any positive non-constant solution to~\eqref{EL0} gives rise to an optimal function for~\eqref{mu-Q}. We recall that $\overline{\mu}(p,\lambda)=\mu(p,1,\lambda)$. Since Proposition~\ref{Prop0:Reparametrization} applies to $\theta=1$, in that case an optimal function solves~\eqref{EL0}, up to a multiplicative constant. Let us assume that there is some continuous branch $\lambda\mapsto u_\lambda$ of optimal functions in the case $\theta=1$, which solve~\eqref{EL0}. Assume moreover that there is no other non-constant solution of~\eqref{EL0}. Then
\[
\overline{\mu}(p,\lambda)=(p-2)\,\nrmS{u_\lambda}p^{p-2}
\]
and we can introduce
\[
\nu(p,\lambda):=\frac{\nrmS{\nabla u_\lambda}2^2}{\nrmS{u_\lambda}2^2}\,.
\]
Here we assume that $\lambda$ is in some interval $\mathcal I$. By testing~\eqref{EL0} with $u_\lambda$, we know that
\[
\lambda+(p-2)\,\nu(p,\lambda)=\overline{\mu}(p,\lambda)\,\frac{\nrmS{u_\lambda}p^2}{\nrmS{u_\lambda}2^2}
\]
and learn that the curve $\Lambda\mapsto\mu(p,\theta,\Lambda)$ is contained in the set
\be{Lambda-M}
\Big\{\big(\Lambda(p,\theta,\lambda),M(p,\theta,\lambda)\big)\,:\,\lambda\in\mathcal I\Big\}\quad\mbox{where}\quad\left\lbrace
\begin{array}{l}
\Lambda(p,\theta,\lambda):=\theta\,\lambda-(1-\theta)\,(p-2)\,\nu(p,\lambda)\,,\\[6pt]
M(p,\theta,\lambda):=\theta\,\overline{\mu}(p,\lambda)^{1/\theta}\(\lambda+(p-2)\,\nu(p,\lambda)\)^{1-1/\theta}\,.
\end{array}
\right.
\ee
\begin{remark} If we consider the bifurcation from constant functions, which occurs for $\theta=1$ at $\lambda=\mu=d$, then by considering the corresponding branch $\lambda\mapsto u_\lambda$ of positive non-constant solutions to~\eqref{EL0}, we may notice that $\lim_{\lambda\to d_+}\overline{\mu}(p,\lambda)=d$, $\lim_{\lambda\to d_+}\nu(p,\lambda)=0$ and, as a consequence
\be{BifurcationPoint}
\lim_{\lambda\to d_+}\Lambda(p,\theta,\lambda)=\lim_{\lambda\to d_+}M(p,\theta,\lambda)=d\,\theta\,.
\ee
The reparametrization~\eqref{Lambda-M} is also consistent with Lemma~\ref{Lem:Asymptotics} as $\lambda\to+\infty$. From the case $\theta=1$, we know that $\overline{\mu}(p,\lambda)\sim\kappa(p,1)\,\lambda^{1-\theta_\star}$ so that $\nu(p,\lambda)\sim c\,\lambda$ for some $c>0$ and recover that $\Lambda(p,\theta,\lambda)\sim\kappa(p,\theta)\,M(p,\theta,\lambda)^{1-\theta_\star/\theta}$ as $\lambda\to+\infty$ using $\big(\kappa(p,\theta)\big)^\theta=\kappa(p,1)\,\theta^\theta/(\theta-\theta_\star)^{\theta-\theta_\star}$.

The case $d=1$ is of particular interest as all positive non-constant solutions of~\eqref{EL0} are known, exist only in $\mathcal I=(d,+\infty)$ and are uniquely for any given $\lambda\in\mathcal I$. However, for a given value of~$\Lambda$,~\eqref{EL} may have several non-constant solutions, even if $d=1$ or up to a rotation in dimension $d\ge1$. Hence the reparametrization in terms of $\lambda$ has to be reserved for numerical computations or a local analysis of the branches and $\mu(p,\theta,\Lambda)$ can be estimated only after optimizing $\mathcal Q_{p,\theta,\lambda}[u]$ on all admissible solutions $u$ of~\eqref{EL0}.
\end{remark}

\subsection{Comparison with inequalities on the Euclidean space}\label{Sec:GN}

We investigate the regime of $\theta$ close to $\theta_\star$ and prove that for $p>2$ but not too large, symmetry breaking occurs for $\lambda<d\,\theta$.
\begin{proposition}\label{Prop:SymmetryBreaking} Let $d\ge1$. There is some $p_\ast\in(2,2^*)$ such that for any $p\in(2,p_\ast)$ and any $\theta\in(\theta_\star,\theta_\star+\varepsilon)$ for some $\varepsilon(p)>0$, we have
\[
\frac{\mu(p,\theta,d\,\theta)}{d\,\theta}<1\,.
\]
\end{proposition}
\begin{proof} Taking the logarithm of both sides in~\eqref{Ineq:Rd}, we have
\[
\frac{\theta_\star}2\,\log\(\nrm{\nabla f}{\mathrm L^2(\R^d)}^2\)+\frac{1-\theta_\star}2\,\log\(\nrm f{\mathrm L^2(\R^d)}^2\)\ge\log\(\mathcal C_{\mathrm{GNS}}(p)\)+\frac1p\,\log\(\ird{|f|^p}\)\,.
\]
This inequality for $p>2$ becomes an equality at $p=2$, with $\mathcal C_{\mathrm{GNS}}(2)=1$. Using $\theta_\star'(2)=d/4$ where $'$ denotes the derivative with respect to $p$, we obtain
\[
\nrmR f2^2\,\log\(\frac{\nrmR{\nabla f}2^2}{\nrmR f2^2}\)\ge\frac2d\ird{|f|^2\,\log\(\frac{|f|^2}{\nrmR f2^2}\)}+\frac8d\,\mathcal C_{\mathrm{GNS}}'(2)\,\nrmR f2^2
\]
by taking the derivative with respect to $p$ at $p=2$, which the logarithmic Sobolev in scale invariant form with optimal constant: see for instance Theorem~2 of~\cite{MR479373}, Inequality~(2.3) of~\cite{MR0109101} or Inequality~(26) of~\cite{MR1132315}. Hence we deduce that $\mathcal C_{\mathrm{GNS}}'(2)=\frac d8\,\log(\pi\,d\,e/2)$. Altogether, this means that
\[
\lim_{\lambda\to+\infty}\frac{\mu(p,\theta_\star,\lambda)}{d\,\theta_\star}\sim\frac{2\,\pi\,e}{d\,|\S^d|^{2/d}}=\frac{2\,e}d\(\frac{\Gamma\big(\tfrac{d+1}2\big)}{2\,\sqrt\pi}\)^{2/d}\,.
\]
As a function of $d$, the right-hand side is increasing with limit $1$ as $d\to+\infty$ according to Stirling's formula, which proves that
\be{mu:theta_star}
\lim_{p\to2^+}\frac{\mu(p,\theta_\star,\lambda)}{d\,\theta_\star}<1\quad\forall\,(d,\lambda)\in\N\times(0,+\infty)\,.
\ee

We apply~\eqref{mu:theta_star} with $\lambda=d\,\theta$ and use the continuity of $\mu(p,\theta,\lambda)$ first with respect to~$p$ and then with respect to $\theta$, as a consequence of the reparametrization.\end{proof}

To get an upper estimate of $\mathcal C_{\mathrm{GNS}}(p)$ in an arbitrary dimension $d\ge1$, we can use concentrating Gaussian test functions. Using~\eqref{Bound:thetastar}, we obtain
\be{gpd}
\mu(p,\theta_\star,\lambda)\le\kappa(p,\theta_\star)=\tfrac{p-2}{|\S^d|^{2/d}}\,\mathcal C_{\mathrm{GNS}}(p)^{2/\theta_\star}\le\mathsf g(p,d):=\tfrac{d\,(p-2)\,p^\frac2{p-2}}{2^\frac{(d+2)\,p-4}{d\,(p-2)}\,\pi^\frac1d}\,\Gamma\(\tfrac{d+1}2\)^{2/d}
\ee
for any $\lambda>0$. The condition $\mathsf g(p,d)<d\,\theta_\star$ shows that for $\theta-\theta_\star>0$ small, symmetry breaking occurs for $\lambda<d\,\theta$. Checking the condition $\mathsf g(p,d)<d\,\theta_\star$ is easy to implement numerically. Solving $\mathsf g(p,d)=d\,\theta_\star$ gives a lower estimate on $p_\ast$. See Section~\ref{Sec:numdge3} for more details.

\subsection{Consequences for the main results}\label{Sec:Main2}

With the estimates of Section~\ref{Sec:More}, we can complete the proofs of Section~\ref{Sec:Main1}.

\medskip\noindent{\bf Proof of Theorem~\ref{Thm1}, part (2/2).} For large values of $\lambda$, the estimate of $\mu(p,\theta,\lambda)$ given by~\eqref{asymptotics} is the result of Lemma~\ref{Lem:Asymptotics} with $\kappa=\kappa(p,\theta)$.\hfill\ $\square$

\medskip\noindent{\bf Proof of Theorem~\ref{Thm2}, part (2/2).} It is easy to rephrase the result of Proposition~\ref{Prop:SymmetryBreaking} as the fact that $\lambda_{\mathrm T}(p,\theta)<d\,\theta$ if $\theta-\theta_\star\ge0$ is taken small enough and $p$ is in the  appropriate range.\hfill\ $\square$

\section{A formal Taylor expansion near the bifurcation point}\label{Sec:Taylor}

We consider the minimization problem~\eqref{mu-Q} and investigate as in~\cite{0951-7715-27-3-435} the behaviour of the branch that bifurcates from the constant functions by formally expanding $\mathcal Q_{p,\theta,\lambda}$ on the lowest spherical harmonics in a neighbourhood of $\lambda=d\,\theta$. The reader who is not interested in the discussion for an arbitrary dimension $d\ge1$ will find in Section~\ref{Sec:1bifurcation} a shorter version specific to the case $d=1$.

In order to avoid heavy notations, assume that $p$ and $\theta$ are given, let us write $\Lambda(\lambda)=\Lambda(p,\theta,\lambda)$ and $M(\lambda)=M(p,\theta,\lambda)$ for the reparametrization of~\eqref{Lambda-M}, and consider the parametric curve $\mathcal C:\,\lambda\mapsto\big(\Lambda(\lambda),M(\lambda)\big)$. If $'$ denotes the derivative with respect to $\lambda$, let
\be{slope}
\mathsf s(\lambda):=\frac{M'(\lambda)}{\Lambda'(\lambda)}\,.
\ee
We shall also need
\be{theta0}
\theta_0:=\frac{(d+2)\,(d+3)\,(p-2)}{2\,(p^2+2\,p-6)+d\,(p^2+6\,p-12)-d^2\,(p-2)^2}\,.
\ee
\begin{proposition} Let $d\ge1$, $p\in(2,+\infty)$ if $d=1$, $2$, $p\in(2,2^*]$ if $d\ge3$ and $\theta\ge\theta_\star$. Assuming that $\mathcal C$ is analytic, there is some $\epsilon>0$ such that the parametric curve $\mathcal C:\,[d,d+\epsilon)\to(\R^+)^2$ has the following properties:
\begin{enumerate}
\item If $\theta\neq\theta_0$, the curve $\mathcal C$ bifurcates from $(d\,\theta,d\,\theta)$ tangentially to $\mu=\lambda$, \emph{i.e.}, $M(d)=\Lambda(d)=d\,\theta$ and $\mathsf s(d)=1$,
\item The curve $\mathcal C$ is concave and below the line $\mu=\lambda$ on a right neighbourhood of the bifurcation point, \emph{i.e.}, $\Lambda'(\lambda)>0$, $M'(\lambda)>0$ and $\mathsf s'(\lambda)<0$ for any $\lambda\in(d,d+\epsilon)$ if $\theta>\theta_0$.
\item The curve $\mathcal C$ is convex and above the line $\mu=\lambda$ on a left neighbourhood of the bifurcation point, \emph{i.e.}, $\Lambda'(\lambda)<0$, $M'(\lambda)<0$ and $\mathsf s'(\lambda)>0$ for any $\lambda\in(d,d+\epsilon)$ if $\theta<\theta_0$.\end{enumerate}
\end{proposition}
\begin{proof} The proof relies on a formal Taylor expansion in which we neglect higher order terms without justification, based on our analyticity assumption. Let us start by introducing the lowest spherical harmonic function that are needed to obtain the leading order terms.

On the sphere $\S^d$, we consider the Laplace-Beltrami operator restricted to the functions depending only on the azimuthal angle $\zeta\in[0,\pi]$
\[
\mathcal Lf:=(\sin\zeta)^{1-d}\,\frac d{d\zeta}\((\sin\zeta)^{d-1}\,\frac{df}{d\zeta}\)\,.
\]
On $[0,\pi]$, we consider the probability measure
\[
d\sigma_\zeta=\frac1{Z_d}\,(\sin\zeta)^{d-1}\,d\zeta\quad\mbox{where}\quad Z_d:=\int_0^\pi(\sin\zeta)^{d-1}\,d\zeta=\frac{\sqrt\pi\;\Gamma(\frac d2)}{\Gamma(\frac{d+1}2)}\,.
\]
The spherical harmonic functions
\[
\varphi_0(\zeta):=1\,,\quad\varphi_1(\zeta)=\sqrt{d+1}\,\cos\zeta\,,\quad\varphi_2=\sqrt{\frac{d+3}{2\,d}}\,(d+1)\,\((\cos\zeta)^2-1\)\,,
\]
are normalised eigenfunctions in $\mathrm L^2([0,\pi],d\sigma_\zeta)$ of $\mathcal L$ with eigenvalues $\lambda_0=0$, $\lambda_1=d$, $\lambda_2=2\,d+2$ such that
\begin{align*}
&\isd{|\varphi_1|^4}=3\,\frac{d+1}{d+3}\,,\quad\isd{|\varphi_1|^2\,\varphi_2}=\sqrt{\frac{2\,d}{d+3}}\,,\quad\isd{|\varphi_{i}|^2}=1\quad\forall\,i=0\,,\,1\,,\,2\,,\\
&\isd{\varphi_1}=\isd{\varphi_2}=\isd{\varphi_1\,\varphi_2}=\isd{|\varphi_2|^2\,\varphi_1}=\isd{\varphi_1^3}=0\,.
\end{align*}
According to Section~\ref{Sec:SymmetryBreaking}, the Taylor expansion of $\mathcal Q_{p,\theta,\Lambda}[1+\varepsilon\,\varphi_1]$ at order $\varepsilon^2$ shows that $\mu(p,\theta,\Lambda)<\Lambda$ if $\Lambda>d\,\theta$, while $\mathcal Q_{p,\theta,\Lambda}[u]$ is linearly stable at $u=1$ if $\Lambda<d\,\theta$. Moreover, we know from~\eqref{BifurcationPoint} that $\Lambda(d)=M(d)=d\,\theta$, so that $\mu(p,\theta,d\,\theta)=d\,\theta$. In order to understand the bifurcation of $\Lambda\mapsto\mu(p,\theta,\Lambda)$ at $\Lambda=d\,\theta$, we have to perform a Taylor expansion at higher order. The unique positive solution of~\eqref{EL0} which is constant is $u=\big(\lambda/(p-2)\big)^{1/(p-2)}$ and we consider the expansion with respect to $\lambda$ in a neighbourhood of $\lambda=d$. Let
\be{uvarepsilon}
u_\varepsilon=\big(\tfrac d{p-2}\big)^\frac1{p-2}\(1+\varepsilon\,\mathsf a\,\varphi_1+\varepsilon^2\,\mathsf b\,\varphi_2\)\,.
\ee
A direct computation shows that
\[
\begin{aligned}
\big(\tfrac d{p-2}\big)^{-\frac2{p-2}}\,\nrmS{\nabla u_\varepsilon}2^2&=d\,\mathsf a^2+2\,(d+1)\,\mathsf b^2\,,\\
\big(\tfrac d{p-2}\big)^{-\frac2{p-2}}\,\nrmS{u_\varepsilon}2^2&=1+\mathsf a^2+\mathsf b^2\,,\\
\big(\tfrac d{p-2}\big)^{-\frac2{p-2}}\,\nrmS{u_\varepsilon}p^2&=1+(p-1)\(\mathsf a^2\,\varepsilon^2+\mathsf c\,\varepsilon^4\)\quad\mbox{with}\quad\mathsf c:=\mathsf b^2+(p-2)\,\sqrt{\tfrac{2\,d}{d+3}}\,\mathsf a^2\,\mathsf b-\tfrac{(p-2)\,(d+p)}{2\,(d+3)}\,\mathsf a^4\,.
\end{aligned}
\]

With these preliminaries in hand, our strategy can be decomposed into three steps and goes as follows.\\[-24pt]
\begin{enumerate}
\item If $\theta=1$, we know that the curve $\mathcal C$ bifurcates from $(\lambda,\mu)=(d,d)$ as a concave function $\lambda\mapsto\mu(\lambda)$, to the right and below the line $\mu=\lambda$. We can compute an approximation on the lowest spherical harmonic function by minimizing $\mathcal Q_{p,1,\lambda}$ with respect to the corresponding coefficients.
\item We use the reparametrization of Section~\ref{Sec:Reparametrization} to obtain a curve which bifurcates from $(\lambda,\mu)=(d\,\theta,d\,\theta)$.
\item We discuss the behaviour of $\mathcal C$ in the neighbourhood of the bifurcation point depending on $\theta$.
\end{enumerate}

\smallskip\noindent{\bf\it Step 1: the case $\theta=1$.}~In view of the Taylor expansion of Section~\ref{Sec:SymmetryBreaking}, let us take $\lambda=\lambda(\varepsilon):=d+\varepsilon^2$ and consider $u_\varepsilon$ given by~\eqref{uvarepsilon}. A minimization of $\mathcal Q_{p,1,\lambda(\varepsilon)}[u_\varepsilon]$ with respect to $\mathsf a$ and $\mathsf b$ shows that
\[
\mathcal Q_{p,1,\lambda(\varepsilon)}[u_\varepsilon]=\lambda(\varepsilon)-\frac{(d+2)\,(d+3)\,(p-2)}{2\,d\,(d+1)\,(p-1)\,(2\,d-p\,(d-2))}\,\varepsilon^4+o(\varepsilon^4)
\]
for the optimal choice of $\mathsf a$ and $\mathsf b$ given by
\[
\mathsf a=\sqrt{\frac{(d+2)\,(d+3)}{d\,(d+1)\,(p-1)\,(2\,d-p\,(d-2))}}\quad\mathrm{and}\quad\mathsf b=\frac{\sqrt{d\,(d+3)}}{\sqrt2\,(d+1)\,(2\,d-p\,(d-2))}\,.
\]

\smallskip\noindent{\bf\it Step 2: the reparametrization, the case $\theta\neq1$.}~According to Section~\ref{Sec:Reparametrization} and with the simplified notation $\mu(\lambda)=\overline{\mu}(p,\lambda)$ and $\nu(\lambda)=\nu(p,\lambda)$, we have
\be{eq:14-1}
\nu\big(\lambda(\varepsilon)\big)=\frac{\nrmS{\nabla u_\varepsilon}2^2}{\nrmS{u_\varepsilon}2^2}+o(\varepsilon^4)=\alpha\,\varepsilon^2+\frac12\,\beta\,\varepsilon^4+o(\varepsilon^4)
\ee
where
\[
\alpha=\frac{(d+2)\,(d+3)}{(d+1)\,(p-1)\,(2\,d-p\,(d-2))}\quad\mbox{and}\quad\beta=\frac{(d+3)\,\big(d^2\,(d+1)\,(p-1)^2-(d+2)^2\,(d+3)\big)}{d\,(d+1)^2\,(p-1)^2\,(2\,d-p\,(d-2))^2}\,,
\]
and $\mu\big(\lambda(\varepsilon)\big)=\mathcal Q_{p,1,\lambda(\varepsilon)}[u_\varepsilon]$, that is,
\be{eq:14-1mu}
\mu\big(\lambda(\varepsilon)\big)=d+\varepsilon^2-\gamma\,\varepsilon^4+o(\varepsilon^4)\quad\mbox{where}\quad\gamma=\frac{(d+2)\,(d+3)\,(p-2)}{2\,d\,(d+1)\,(p-1)\,(2\,d-p\,(d-2))}\,.
\ee

\smallskip\noindent{\bf\it Step 3: discussion.}~In our simplified notation, the reparametrization~\eqref{Lambda-M} can be rewritten as
\[
\Lambda(\lambda)=\theta\,\lambda-(1-\theta)\,(p-2)\,\nu(\lambda)\quad\mbox{and}\quad M(\lambda)=\theta\,\mu(\lambda)^{1/\theta}\(\lambda+(p-2)\,\nu(\lambda)\)^{1-1/\theta}\,.
\]
{}From the analyticity assumption, the derivatives can be computed order by order in terms of $\varepsilon^2$ and we obtain that
\[
\Lambda'(d)=\theta-(1-\theta)\,(p-2)\,\alpha>0
\]
if and only if $\theta>\theta_0$. By direct computation, we have
\[
M'(\lambda)=\mu'(\lambda)\(\frac{\lambda+(p-2)\,\nu(\lambda)}{\mu(\lambda)}\)^{1-1/\theta}+(\theta-1)\(\frac{\mu(\lambda)}{\lambda+(p-2)\,\nu(\lambda)}\)^{1/\theta}\,\big(1+(p-2)\,\nu'(\lambda)\big)
\]
and for $\lambda=d$, $\mu(d)=d$, while we learn that $\nu(d)=0$ and $\mu'(d)=1$ from~\eqref{eq:14-1} and~\eqref{eq:14-1mu} respectively. Altogether, we obtain
\[
M'(d)=\theta-(1-\theta)\,(p-2)\,\alpha=\Lambda'(d)\,.
\]
As a consequence, the slope $\mathsf s$ of the curve $\mathcal C$ defined by~\eqref{slope} satisfies
\[
\mathsf s(d)=1\quad\mbox{if}\quad\theta\neq\theta_0\,.
\]
If $\theta\neq\theta_0$, let us observe that
\[
\mathsf s'(d)=\frac{\mathsf s'(d)}{\mathsf s(d)}=\frac{M''(d)}{M'(d)}-\frac{\Lambda''(d)}{\Lambda'(d)}=\frac{M''(d)-\Lambda''(d)}{\Lambda'(d)}
\]
where
\[
\Lambda''(d)=-(1-\theta)\,(p-2)\,\beta\quad\mbox{and}\quad M''(d)=-\mu''(d)+(p-2)\,\frac{1-\theta}{d\,\theta}\((p-2)\,\alpha^2-d\,\theta\,\beta\)\,.
\]
Altogether, we have the identity
\[
\mathsf s'(d)\,\Lambda'(d)=-\,(p-2)^2\,\frac{1-\theta}{d\,\theta}\,\alpha^2
\]
from which we deduce the concavity (resp.~convexity) property of $\mathcal C$ if $\theta>\theta_0$ (resp.~$\theta<\theta_0$).
\end{proof}

\section{The one-dimensional case}\label{Sec:d=1}

\subsection{Computation of the branch with \texorpdfstring{$\theta=1$}{theta=1}}\label{Sec:1dtheta=1}

Let us consider the solution of~\eqref{EL0}, up to a multiplication by $\big(\lambda/(p-2)\big)^\frac1{p-2}$. The $2\,\pi$-periodic solution of
\be{Eqn:lambda}
-\,u''+\frac\lambda{p-2}\(u-u^{p-1}\)=0\quad\mbox{on}\quad\R
\ee
is changed into the $T$-periodic solution of
\be{Eqn:v}
-\,v''+v=v^{p-1}\quad\mbox{on}\quad\R
\ee
with $T=2\,\pi\,\sqrt{\lambda/(p-2)}$ if $u(x)=v\big(\sqrt{\lambda/(p-2)}\,x\big)$. Since
\[
E=\frac12\,|v'|^2-\Pot(v)\quad\mbox{where}\quad\Pot(v):=\tfrac12\,|v|^2-\tfrac1p\,|v|^p
\]
is constant, we learn that
\[
|v'|={\sqrt{2\(E+\Pot(v)\)}}
\]
and can compute the period of a periodic positive function as
\[
T(E)=2\int_0^{T(E)/2}dx=2\iE{\frac{dw}{\sqrt{2\(E+\Pot(w)\)}}}
\]
where
\[
I(E):=\big\{w>0\,:\,E+\Pot(w)>0\big\}\quad\forall\,E\in\(-\,\theta_\star,0\)
\]
and $\theta_\star=(p-2)/(2\,p)=\min_{w>0}\Pot(w)$. Changing variables with $\lambda=\overline{\lambda}(E)$, we can compute
\begin{align}
&\is{|u|^q}=\frac{\sqrt{p-2}}{2\,\pi\,\sqrt\lambda}\int_0^{T(E)}|v|^q\,dx=\frac{\sqrt{p-2}}{\pi\,\sqrt\lambda}\,J_q(E)\,,\label{Jq}\\
&\is{|u'|^2}=\frac{\sqrt\lambda}{2\,\pi\,\sqrt{p-2}}\int_0^{T(E)}|v'|^2\,dx=\frac{\sqrt\lambda}{\pi\,\sqrt{p-2}}\,K(E)\,,\label{K}
\end{align}
with
\[
J_q(E):=\iE{\frac{w^q}{\sqrt{2\(E+\Pot(w)\)}}}\,dw\quad\mbox{and}\quad K(E):=\iE{\sqrt{2\(E+\Pot(w)\)}}\,dw\,.
\]
According to~\cite{dolbeault2024monotonicity} (also see earlier references therein), the function $T:[-\,\theta_\star,0)\to\R^+$ is monotone increasing with
\[
\lim_{E\to(-\,\theta_\star)_+}T(E)=2\,\pi\quad\mbox{and}\quad\lim_{E\to0_-}T(E)=+\infty\,.
\]
The concave curve $\lambda\mapsto\overline{\mu}(p,\lambda)$ can be parametrized by $E$. With a slight abuse of notations, we can consider $\overline{\lambda}$ and $\overline{\mu}$ as functions of $E$ and write
\begin{align*}
&\overline{\lambda}(E)=(p-2)\(\frac{T(E)}{2\,\pi}\)^2\,,\\
&\overline{\mu}(E)=\frac{(p-2)\is{|u'|^2}+\overline{\lambda}(E)\is{|u|^2}}{\(\is{|u|^p}\)^{2/p}}=\overline{\lambda}(E)\(\frac{\sqrt{p-2}}{\pi\,\sqrt{\overline{\lambda}(E)}}\)^{1-\frac2p}
\frac{K(E)+J_2(E)}{J_p(E)^{2/p}}\,,
\end{align*}
For later purpose, it is also convenient to define
\[
\overline{\nu}(E):=\frac{\nrms{u'}2^2}{\nrms u2^2}=\frac{\overline{\lambda}(E)}{p-2}\,\frac{K(E)}{J_2(E)}\,.
\]
Altogether, we have
\be{Parametrization:E}
\begin{aligned}
&\overline{\lambda}(E)=(p-2)\(\frac{T(E)}{2\,\pi}\)^2\,,\\
&\overline{\mu}(E)=(p-2)\(\frac{T(E)}{2\,\pi}\)^2\,\(\frac{T(E)}2\)^{-\frac{p-2}p}\,\frac{K(E)+J_2(E)}{J_p(E)^{2/p}}\,,\\
&\overline{\nu}(E)=\(\frac{T(E)}{2\,\pi}\)^2\,\frac{K(E)}{J_2(E)}\,.
\end{aligned}
\ee
\begin{remark} Notice that
\be{KJJ}
\overline{\mu}(E)=\big(\overline{\lambda}(E)\big)^{1-\theta_\star}\(\frac{p-2}{\pi^2}\)^{\theta_\star}\,\frac{K(E)+J_2(E)}{J_p(E)^{2/p}}\,.
\ee
Since $\lim_{E\to0_-}T(E)=+\infty$ and $J_2(E)$, $J_p(E)$ and $K(E)$ converge as $E\to0_+$ respectively to
\be{KJJid}
J_2(0)=\frac14\,\frac{p+2}{p-2}\(\frac p2\)^\frac 2{p-2}\frac{\sqrt\pi\;\Gamma\big(\tfrac p{p-2}\big)}{\Gamma\big(\tfrac p{p-2}+\tfrac12\big)}\,,\quad J_p(0)=\frac{2\,p}{p+2}\,J_2(0)\,,\quad K(0)=\frac{p-2}{p+2}\,J_2(0)\,,
\ee
we recover the result of Lemma~\ref{Lem:Asymptotics} in dimension $d=1$ if $\theta=1$. For details on this computation, see for instance the proof of Lemma~3 in~\cite{delPino20102045}.
\end{remark}

\subsection{Computation of the branch with \texorpdfstring{$\theta<1$}{theta<1}: reparametrization}\label{Sec:Reparam}

If $\theta\in(\theta_\star,1)$, we know from~Section~\ref{Sec:Reparametrization} that any optimal function for~\eqref{mu-Q} is, up to a multiplication by a constant, a solution of~\eqref{Eqn:lambda} with
\[
\lambda=\frac1\theta\(\Lambda+(1-\theta)\,(p-2)\,\frac{\nrms{u'}2^2}{\nrms u2^2}\)
\]
which provides us with the reparametrization $\big(\Lambda,\mu(p,\theta,\Lambda)\big)$ by $E\mapsto\big(\lambda(\theta,E),\mu(\theta,E)\big)$
\be{Reparametrisation}
\hspace*{-8pt}\begin{aligned}
&\lambda(\theta,E):=\theta\,\overline{\lambda}(E)-(1-\theta)\,(p-2)\,\overline{\nu}(E)=\(\theta-(1-\theta)\,\frac{K(E)}{J_2(E)}\)\overline{\lambda}(E)\,,\\
&\mu(\theta,E):=\big(\lambda(\theta,E)+(p-2)\,\overline{\nu}(E)\big)\(\frac{\nrms{u_E}2^2}{\nrms{u_E}p^2}\)^\frac1\theta=\(\lambda(\theta,E)+\overline{\lambda}(E)\,\frac{K(E)}{J_2(E)}\)\,\(\frac{\sqrt{p-2}}{\pi\,\sqrt{\overline{\lambda}(E)}}\)^\frac{p-2}{p\,\theta}\!\(\frac{J_2(E)}{J_p(E)^{2/p}}\)^\frac1\theta\,,
\end{aligned}\hspace*{-8pt}
\ee
where $u_E$ denotes the solution of~\eqref{Eqn:lambda} corresponding to $\lambda=\overline{\lambda}(E)$ and $J_2$ and $K$ are defined respectively in~\eqref{Jq} and~\eqref{K}. Using the computations of Section~\ref{Sec:1dtheta=1}, we obtain
\be{Reparametrisation:E-theta}
\begin{aligned}
&\lambda(\theta,E)=(p-2)\(\theta-(1-\theta)\,\frac{K(E)}{J_2(E)}\)\(\frac{T(E)}{2\,\pi}\)^2\,,\\
&\mu(\theta,E)=\theta\,(p-2)\(1+\frac{K(E)}{J_2(E)}\)\(\frac{T(E)}{2\,\pi}\)^2\,\(\frac{T(E)}2\)^{-\frac{p-2}{p\,\theta}}\(\frac{J_2(E)}{J_p(E)^{2/p}}\)^\frac1\theta\,.
\end{aligned}
\ee
\begin{remark}\label{Rem:mu} As a limit case corresponding to $\theta=\theta_\star$, we read from~\eqref{Reparametrisation:E-theta} that
\[
\mu(\theta_\star,E)=\theta_\star\,(p-2)\,\frac{K(E)+J_2(E)}{J_2(E)}\(\frac{J_2(E)}{J_p(E)^{2/p}}\)^\frac1{\theta_\star}\,.
\]
Using~\eqref{KJJid}, we obtain
\[
\kappa(p,\theta_\star)=\lim_{E\to0_-}\mu(\theta_\star,E)=\frac{(p-2)^2}{p+2}\(\frac{p+2}{2\,p}\)^{\frac4{p-2}}\,\frac{J_2(0)^2}{\pi^2}
\]
with $J_2(0)$ given by~\eqref{KJJ}. Taking into account~\eqref{Bound:thetastar}, this allows us to check that 
\be{CGNS}
\(\mathcal C_{\mathrm{GNS}}(p)\)^{2/\theta_\star}=\(\tfrac{p+2}4\)^\frac{p+2}{p-2}\!\frac{\pi\;\Gamma\(\frac p{p-2}\)^2}{(p-2)\;\Gamma\(\frac p{p-2}+\frac12\)^2}\,.
\ee
By testing~\eqref{Eqn:v} with $v$, we learn that $\int_0^{T(E)}\(|v'|^2+|v|^2-|v|^p\)\,dx=0$ and
\[
K(E)+J_2(E)-J_p(E)=0
\]
using the change of variables $x\mapsto w=v(x)$ as in Section~\ref{Sec:1dtheta=1}. Similarly an integration on $(0,T(E))$ of the identity $0=2\(E+\Pot(v)\)-|v'|^2$ gives the relation
\[
E\,T(E)-K(E)+J_2(E)-\tfrac2p\,J_p(E)=0\,.
\]
As a consequence, we learn that $J_p(E)=K(E)+J_2(E)$,
\[
T(E)=\frac{(p-2)\,J_2(E)-(p+2)\,K(E)}{p\,|E|}=2\,\frac{J_2(E)}{|E|}\,\(\theta_\star-(1-\theta_\star)\,\frac{K(E)}{J_2(E)}\)
\]
so that $\lambda(\theta_\star,E)\ge0$ and $\lambda(\theta,E)=\frac1{1-\theta_\star}\,\big((p-2)\,(\theta-\theta_\star)+(1-\theta)\,\lambda(\theta_\star,E)\big)\,\big(\frac{T(E)}{2\,\pi}\big)^2\ge(p-2)\,\frac{\theta-\theta_\star}{1-\theta_\star}$ using $T(E)\ge2\,\pi$, and
\[
\lambda(\theta_\star,E)=\frac{(p-2)\,\big((p+2)\,K(E)-(p-2)\,J_2(E)\big)^3}{8\,p^3\,\pi^2\,J_2(E)\,E^2}\,.
\]
A lengthy computation in the spirit of the method of~Section~4 of~\cite{dolbeault2024monotonicity} shows that $K'(E)$ and $J_2'(E)$ admit finite limits as $E\to0_-$ so that $(p+2)\,K(E)-(p-2)\,J_2(E)=O(E)$. As a consequence we know that
\[
\lim_{E\to 0_-}\lambda(\theta_\star,E)=0\,.
\]
\end{remark}

\subsection{Comparison with the Gagliardo-Nirenberg constant}\label{Sec:GN1d}

In the $d=1$ case, it is known that equality in the \emph{Euclidean Gagliardo-Nirenberg-Sobolev inequality}~\eqref{Ineq:Rd} is achieved by the function $u(x)=(\cosh x)^{-2/(p-2)}$, $x\in\R$. See for instance~\cite{Dolbeault06082014} and earlier references therein. As a consequence, we recover~\eqref{CGNS}. 
We know from~\eqref{Bound:thetastar} that $\mu(p,\theta_\star,\lambda)\le\kappa(p,\theta_\star)$. Using the concavity of $\lambda\mapsto\mu(p,\theta,\lambda)$ and its continuity with respect to $\theta$ for $\theta-\theta_\star>0$, small, we conclude that symmetry breaking occurs for $\lambda<d\,\theta$ with $d=1$ if
\[
\kappa(p,\theta_\star)=\frac{p-2}{4\,\pi^2}\,\mathcal C_{\mathrm{GNS}}(p)^{2/\theta_\star}<\theta_\star=\frac{p-2}{2\,p}\,.
\]
One can show that this occurs if $p\in(2,p_\star)$ for some $p_\star>2$. See Section~\ref{Sec:critical-theta} for a numerical result.

\section{Some numerical results}\label{Sec:Numerics}

In this section, we collect various plots which illustrate our results. Except in Sections~\ref{Sec:Numerics-Range} and~\ref{Sec:numdge3}, we focus on the case $d=1$ and use the representation of Section~\ref{Sec:1dtheta=1} of the solutions by elliptic integrals. In higher dimensions, more robust numerical methods would be needed for detailed results, which are out of the scope of the present paper, but at least we present evidences that some properties are qualitatively similar to the case $d=1$.

\subsection{Inequalities, symmetry and phase transition}\label{Sec:Numerics-Range}

\begin{figure}[ht]
\begin{center}
\includegraphics[width=4.5cm]{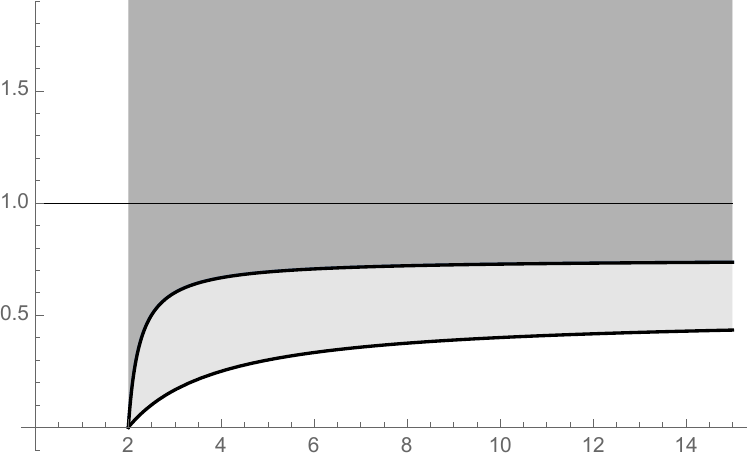}\hspace*{0.4cm}\includegraphics[width=4.5cm]{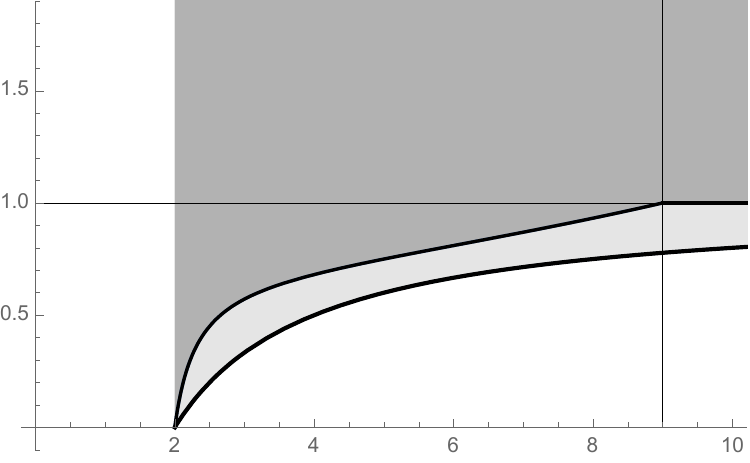}\hspace*{0.4cm}\includegraphics[width=4.5cm]{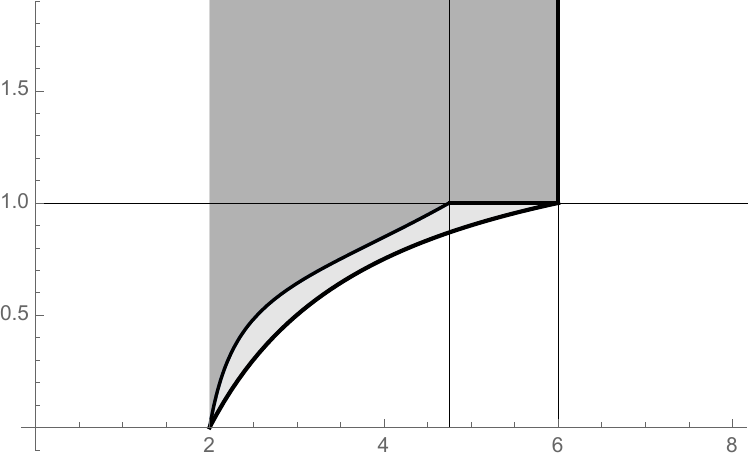}\hspace*{0.4cm}
\caption{\small\label{Fig1} The range of the admissible parameters $\theta$ as a function of $p$ in (from left to right) dimensions $d=1$, $d=2$ and $d=3$ corresponds to the gray area $\theta\ge\theta_\star$ for any $p>2$ if $d=1$ or $d=2$ and $p\in(2,2^*]$ is $d\ge3$. The dark gray subdomain is the region $\theta\ge\theta^\#$ with $p\in(2,+\infty)$ if $d=1$, $2<p\le2^\#$ if $d\ge2$ or $\theta\ge1$ if $d\ge2$ with $2^\#\le p<\infty$ if $d=2$ and $2^\#\le p\le2^*$ if $d\ge3$. In dimensions $d\ge3$, the admissible range is limited to the interval $(2,2^*]$ with $2^*\to2$ as $d\to+\infty$ but otherwise the picture is qualitatively the same as in the case $d=3$.}
\end{center}
\vspace*{-0.5cm}
\end{figure}

The family of inequalities~\eqref{Ineq:GNS2} holds for any $u\in\mathrm H^1(\S^d,d\sigma)$ and $\lambda>0$, for some positive optimal constant $\mu(p,\theta,\lambda)$, if $\theta\ge\theta_\star=d\,(p-2)/(2\,p)$, $p\in(2,+\infty)$ if $d=1$ or $d=2$, and $p\in(2,2^*]$ with $2^*=2\,d/(d-2)$ if $d\ge3$. We have the estimates $\mu(p,\theta,\lambda)\le\lambda$ for any $\lambda>0$, $\mu(p,\theta,\lambda)<\lambda$ if $\lambda>d\,\theta$ and $\mu(p,\theta,\lambda)=\lambda$ if $\lambda\le d\,\theta$, $\theta\ge\theta^\#$ and $p\in(2,2^\#)$, where $2^\#$ is the \emph{Bakry-Emery exponent} and~$\theta^\#$ is the corresponding exponent. If $d=1$, then $2^\#=+\infty$ and $\theta^\#=3\,(p-2)/(4\,p-7)$. If $d\ge2$, we also have $\mu(2^\#,\theta,\lambda)=\lambda$ if $\lambda\le d\,\theta$ and $\theta\ge\theta^\#$. If $p\in(2^\#,2^*$ and $d\ge3$, or $p\in(2^\#,+\infty)$ if $d=2$, we also have $\mu(p,\theta,\lambda)=\lambda$ if $\lambda\le d\,\theta$, $\theta\ge1$. See Fig.~\ref{Fig1}.

The function $\lambda\mapsto\mu(p,\theta,\lambda)$ is monotone non-decreasing and concave. We distinguish two regimes:\\
$\bullet$ \emph{Symmetry:} $\mu(p,\theta,\lambda)=\lambda$, the equality case in~\eqref{Ineq:GNS2} is achieved only by constant functions,\\
$\bullet$ \emph{Symmetry breaking:} $\mu(p,\theta,\lambda)<\lambda$, the equality case in~\eqref{Ineq:GNS2} is not achieved by constant functions.\\
Our goal is to understand whether the \emph{phase transition} from symmetry to symmetry breaking occurs at $\lambda=d\,\theta$ and determines a \emph{second order phase transition}, or for some $\lambda<d\,\theta$, which corresponds to a \emph{first order phase transition}. Notice that a first order phase transition can occur also at $\lambda=d\,\theta$. In Fig.~\ref{Fig1}, dark grey regions correspond to $\theta\ge\min\big\{1,\theta^\#\big\}$ in which we have a second order phase transition.

\subsection{The branch of positive non-constant solutions for \texorpdfstring{$d=1$}{d=1}}

\begin{figure}[ht]
\begin{center}
\includegraphics[width=3.5cm]{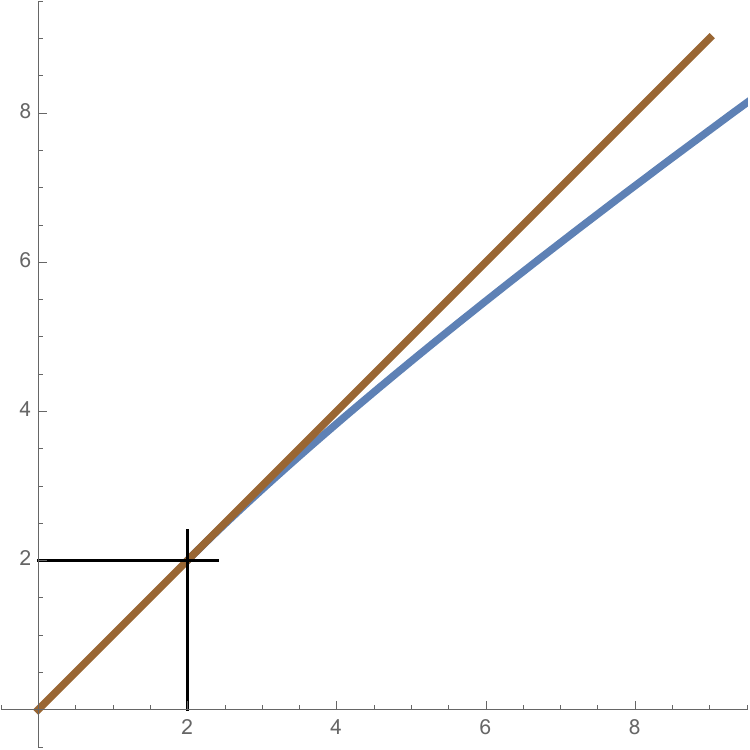}\hspace*{0.3cm}\includegraphics[width=3.5cm]{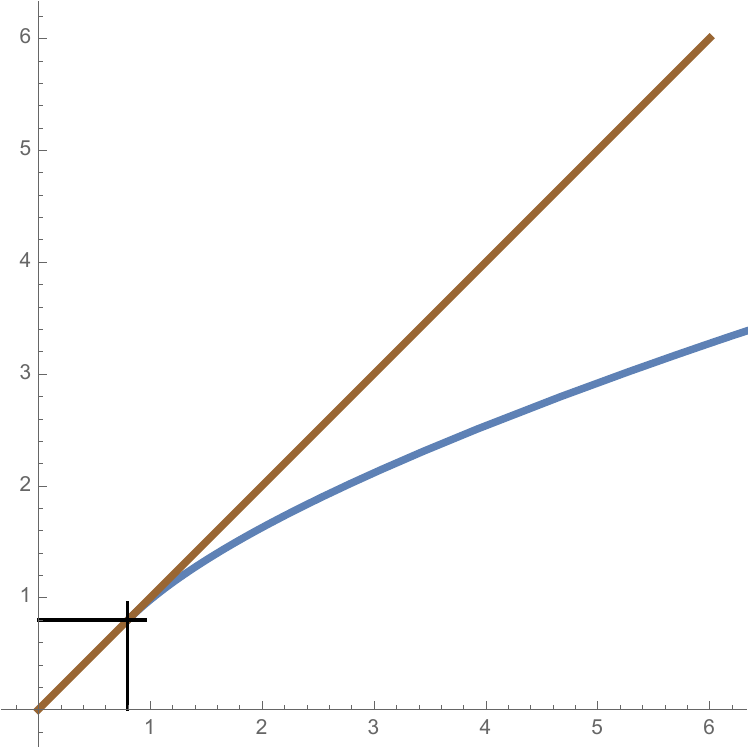}\hspace*{0.3cm}\includegraphics[width=3.5cm]{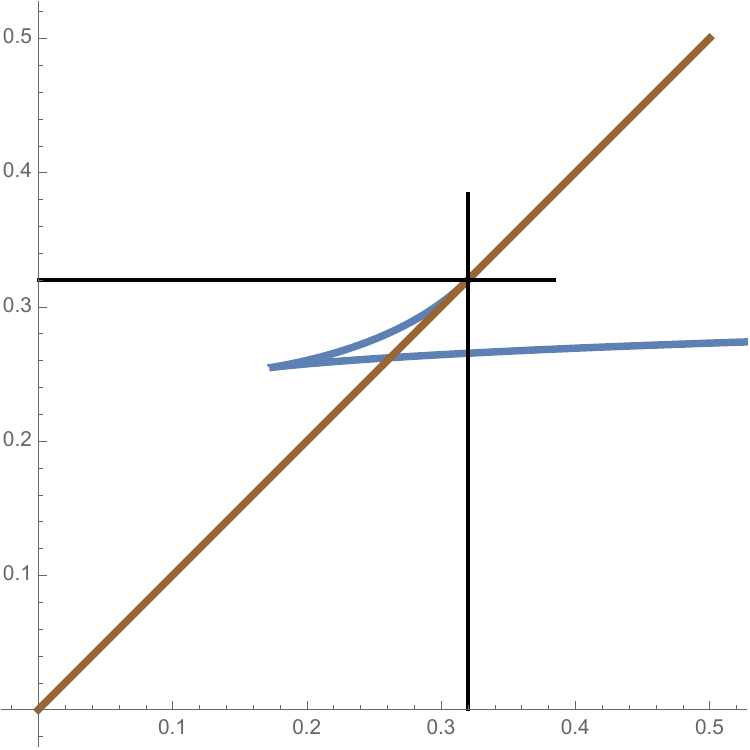}\hspace*{0.3cm}\includegraphics[width=3.5cm]{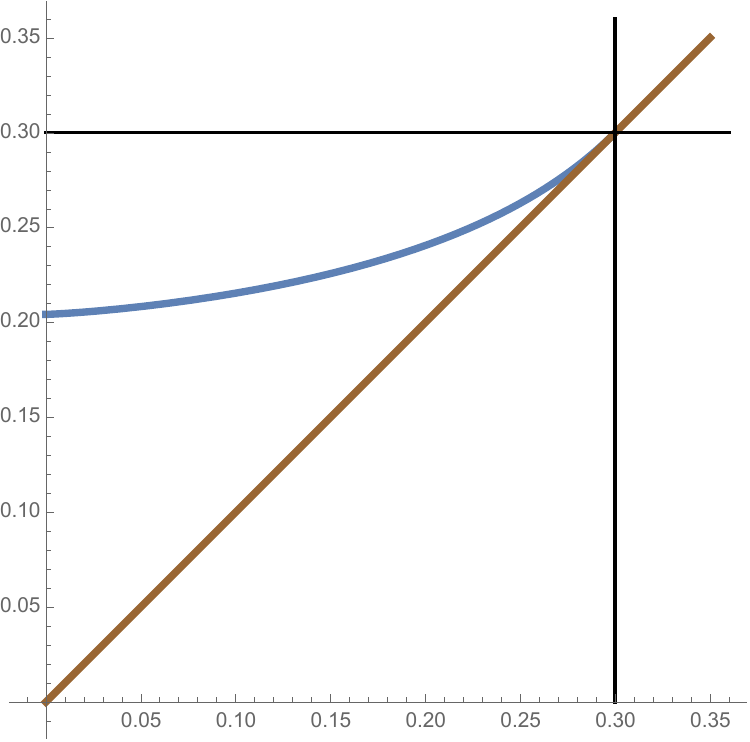}
\caption{\small\label{Fig2} In dimension $d=1$, with $p=5$, the set of non-constant positive solutions is represented in $(\lambda,\mu)$ coordinates (blue curve) for $\theta=2$, $\theta=0.8$, $\theta=0.32$ and $\theta=\theta_\star=0.3$, from left to right, while constant solutions are represented by the brown half-line $\mu=\lambda$. The bifurcation point is at $\lambda=\mu=d\,\theta$ (with $d=1$) in all cases. In the case $\theta=\theta_\star$, we have already seen in Remark~\ref{Rem:mu} that $\lim_{E\to 0_-}\lambda(\theta_\star,E)=0$ and $\lim_{E\to 0_-}\mu(\theta_\star,E)=\kappa(p,\theta_\star)>0$.}
\end{center}
\vspace*{-0.5cm}
\end{figure}

In dimension $d=1$, obtaining all positive non-constant solutions to~\eqref{Eqn:lambda} can be done by solving this Hamiltonian ODE using elliptic integrals involving the Hamiltonian energy~$E$ as a parameter, which parametrizes monotonically $\lambda=\overline{\lambda}(E)$ while $\overline{\mu}(E)=\mu\big(p,1,\lambda(\theta,E)\big)$ is obtained by computing $\overline{\mu}(E)=\overline{\lambda}(E)\,\nrms{u_E}p^{p-2}$ where $u_E$ denotes the $2\,\pi$-periodic solution of energy $E$. The curve $E\mapsto\big(\overline{\lambda}(E),\overline{\mu}(E)\big)$ given by~\eqref{Parametrization:E} represents all solutions corresponding to the Euler-Lagrange equations associated with~\eqref{Ineq:GNS2} for $\theta=1$. If $d=1$ but $\theta\neq1$, all solutions corresponding to the Euler-Lagrange equations associated with~\eqref{Ineq:GNS2} can be represented by a curve $E\mapsto\big(\lambda(\theta,E),\mu(\theta,E)\big)$ given by~\eqref{Reparametrisation:E-theta}. Typical patterns are shown in Fig.~\ref{Fig2}. Depending on the values of $\theta$, the bifurcating curve goes either to the right (higher values of $\lambda$) or to the right. This first qualitative property is explained next.

Coming back to~\eqref{Ineq:GNS1} and~\eqref{Ineq:GNS2}, we may notice that $\overline{\mu}\big(p,\lambda(1,E)\big)=\mu(1,E)$ for any $E>0$ and
\[
\mu\big(p,\theta,\Lambda\big)=\min\Big\{\mu_\ast(\theta,E,\Lambda)\big\}\,:\,E>0\Big\}\quad\mbox{where}\quad\mu_\ast(\theta,E,\Lambda):=\min\{\Lambda,\mu(\theta,E)\}\,.
\]

\subsection{The first bifurcation from constant solutions for \texorpdfstring{$d=1$}{d=1}}\label{Sec:1bifurcation}

\begin{figure}[ht]
\vspace*{-0.25cm}
\begin{center}
\includegraphics[width=6cm]{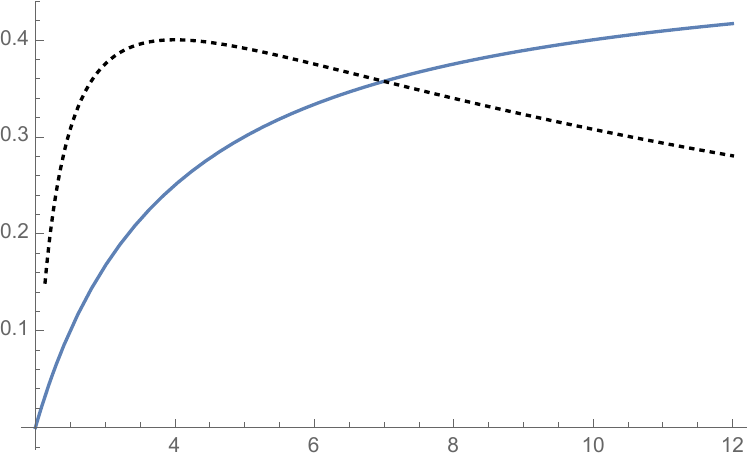}
\caption{\small\label{Fig:Turning} In dimension $d=1$, plot of $p\mapsto\theta_\star(p)$ (blue) and $p\mapsto\theta_0(p)$ (black, dotted).}
\end{center}
\vspace*{-0.5cm}
\end{figure}

The deficit associated with~\eqref{Ineq:GNS2} is linearly stable if $\lambda<d\,\theta$, and linearly unstable if $\lambda>d\,\theta$. As a consequence, a branch of non-constant solutions of the Euler equation bifurcates at $\lambda=d\,\theta$. The case of dimension $d=1$ enters in the general setting of Section~\ref{Sec:Taylor}, but one can give a simpler analysis in that case. Let us start with the case $\theta=1$. The simple \emph{ansatz}
\[
u_\varepsilon(\zeta)=1+\mathsf a\,\varepsilon\,\sqrt2\,\cos \zeta+\mathsf b\,\varepsilon^2\,\sqrt2\,\cos(2\,\zeta)\quad\forall\,\zeta\in[0,2\,\pi)\approx\S^1
\]
allows us to describe the branch of positive non-constant solutions to~\eqref{EL0} near the bifurcation point, at leading order for $\varepsilon>0$ small, as a parametric curve in $(\lambda,\mu)$ coordinates parametrized by $\varepsilon>0$ such that
\[
\lambda(\varepsilon)=1+\varepsilon^2\,.
\]
This \emph{ansatz} has been justified in Section~\ref{Sec:Taylor}. Here we only collect the observations needed for numerical computations in dimension $d=1$. At formal level, we do the expansion of
\[
\mu(\varepsilon):=\frac{(p-2)\,\nrms{u_\varepsilon'}2^2+\lambda(\varepsilon)\,\nrms{u_\varepsilon}2^2}{\nrms{u_\varepsilon}p^2}
\]
and simply observe that $\lambda(\varepsilon)-\mu(\varepsilon)$ achieves its maximum for
\[
\mathsf a=\frac{6}{(p-1)\,(p+2)}\quad\mbox{and}\quad\mathsf b=\frac1{\sqrt2\,(p+2)}\,,
\]
in which case we obtain
\[
\mu(\varepsilon)=1+\varepsilon^2-\gamma\,\varepsilon^4+o\big(\varepsilon^4\big)\quad\mbox{with}\quad\gamma=\frac{3\,(p-2)}{(p-1)\,(p+2)}\,.
\]
Notice that we can eliminate $\varepsilon^2=\lambda-1$ and find that $\overline{\mu}(p,\lambda)=\lambda-\gamma\,(\lambda-1)^2+o\big((\lambda-1)^2\big)$, \emph{i.e.}, the parabola which interpolates $\lambda\mapsto\overline{\mu}(p,\lambda)$ at the bifurcation point. With a similar computation, we can compute
\[
\nu(\varepsilon):=\frac{\nrms{u_\varepsilon'}2^2}{\nrms{u_\varepsilon}2^2}=\alpha\,\varepsilon^2+\beta\,\varepsilon^4+o\big(\varepsilon^4\big)\quad\mbox{where}\quad\alpha=\frac6{(p-1)\,(p+2)}\quad\mbox{and}\quad\beta=\frac{2\,(p^2-2\,p-17)}{(p-1)^2\,(p+2)^2}\,.
\]
So far, these computations are done in the case $\theta=1$. If $\theta\neq1$, we can use them to get an approximation of $\lambda(\theta,E)$ and $\mu(\theta,E)$ using~\eqref{Reparametrisation:E-theta}. For $\varepsilon>0$ small, we obtain the reparametrization
\[
\begin{aligned}
&\lambda=\theta\,(1+\varepsilon^2)-(1-\theta)\,(p-2)\(\alpha\,\varepsilon^2+\beta\,\varepsilon^4\)+o\big(\varepsilon^4\big)\,,\\
&\mu=\theta\((p-2)\(\alpha\,\varepsilon^2+\beta\,\varepsilon^4\)+1+\varepsilon^2\)^{1-\frac1\theta}\(1+\varepsilon^2-\gamma(p)\,\varepsilon^4
\)^\frac1\theta+o\big(\varepsilon^4\big)\,.
\end{aligned}
\]
Expanding again around $\varepsilon=0$, we notice that
\[
\begin{aligned}
&\lambda=\theta+A(\theta,p)\,\varepsilon^2+B(\theta,p)\,\varepsilon^4+o\big(\varepsilon^4\big)\\
&\mu=\theta+A(\theta,p)\,\varepsilon^2+C(\theta,p)\,\varepsilon^4+o\big(\varepsilon^4\big)
\end{aligned}\quad\mbox{with}\quad\left\{
\begin{aligned}
&A(\theta,p)=\theta-(1-\theta)\,(p-2)\,\alpha=\theta-(1-\theta)\,\tfrac{6\,(p-2)}{(p-1)\,(p+2)}\,,\\
&B(\theta,p)=-\,(1-\theta)\,(p-2)\,\beta\,,\\
&C(\theta,p)=\tfrac{1-\theta}{2\,\theta}\,(p-2)^2\,\alpha^2-(1-\theta)\,(p-2)\,\beta-\gamma\,.
\end{aligned}\right.
\]
To understand the behaviour of the branch in the $(\lambda,\mu)$ representation, it is interesting to compute the sign of $A(\theta,p)$. An elementary computation shows that $A(\theta,p)$ is negative if and only if
\[
\theta<\theta_0(p):=\frac{6\,(p-2)}{p^2+7\,p-14}<\theta_\star\,.
\]
This is consistent with~\eqref{theta0}. As in Section~\ref{Sec:Taylor}, we obtain that $\mathsf s'(\lambda)<0$ for $\lambda>d$, close enough to $d$ if $\theta<\theta_0$, so that bifurcations always take place \emph{to the left} in that case. This happens if $2<p<p_0=7$. See~Fig.~\ref{Fig:Turning}.

\begin{figure}[ht]
\vspace*{-0.25cm}
\begin{center}
\includegraphics[height=4cm]{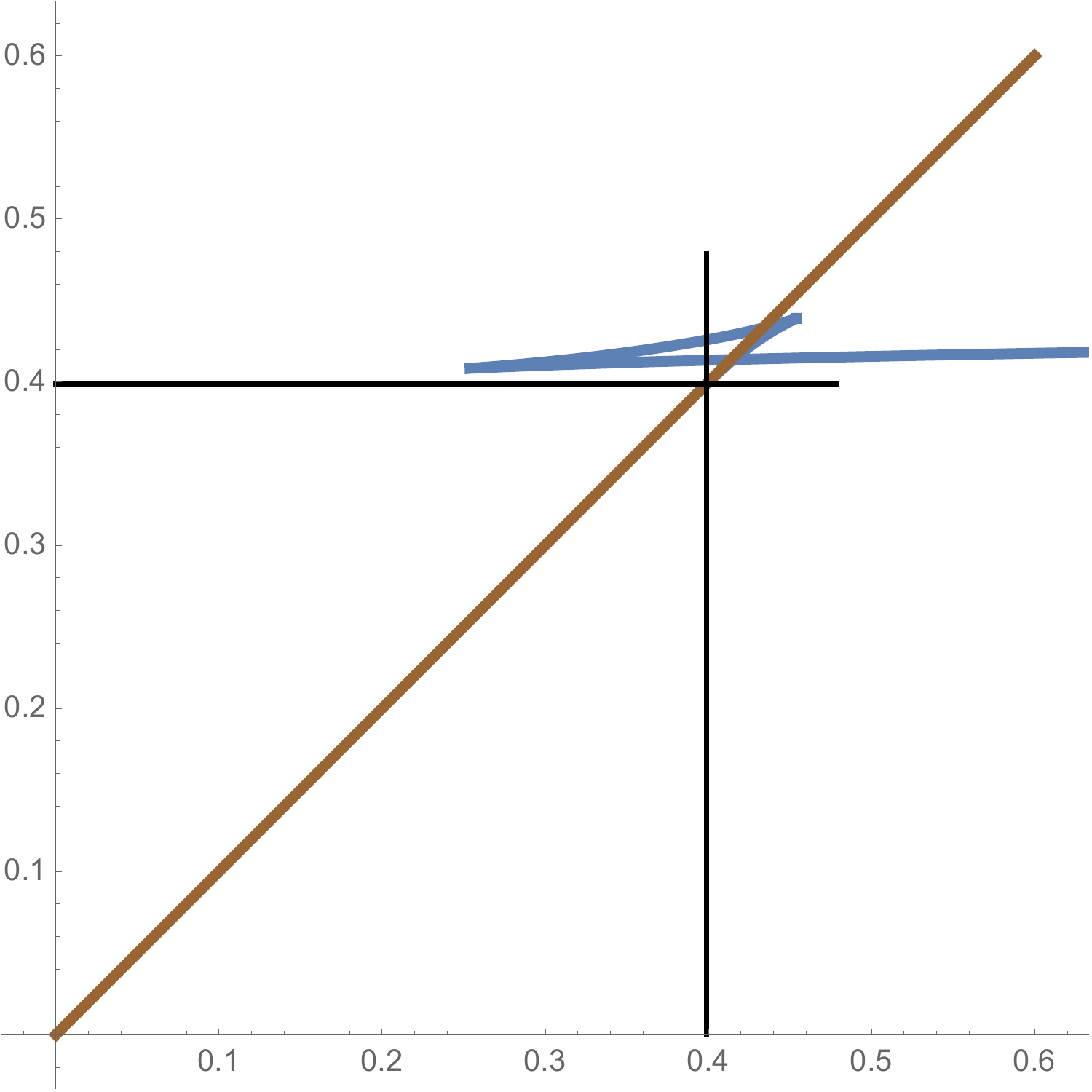}\hspace*{0.5cm}\includegraphics[height=2.8cm]{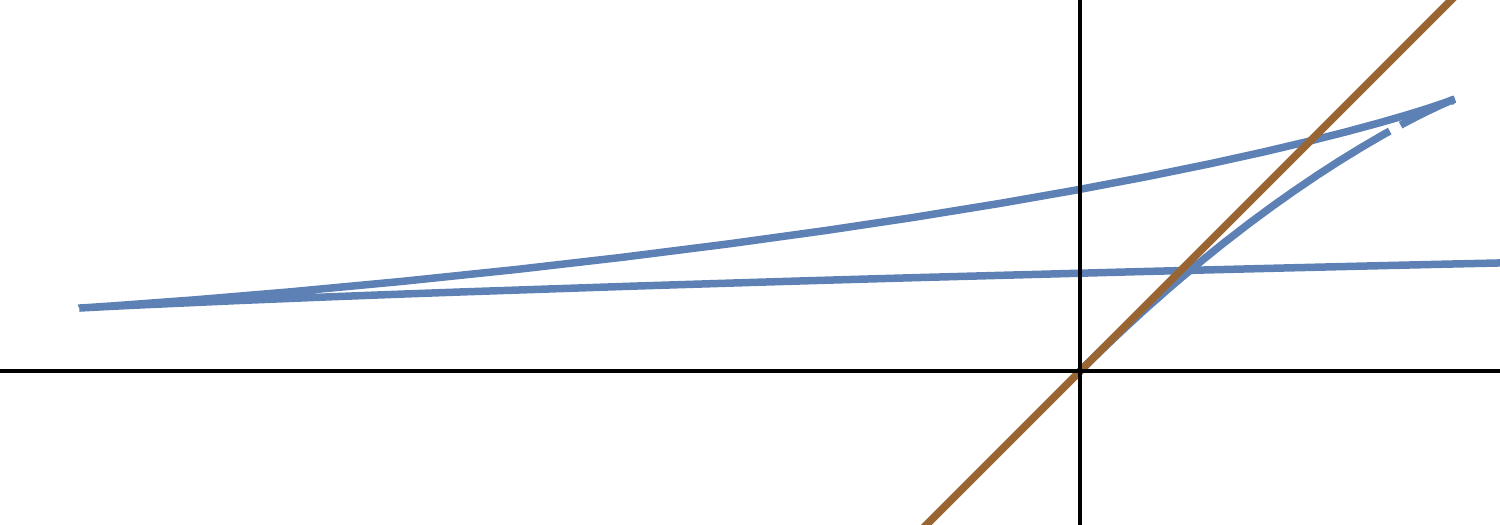}
\caption{\small\label{Fig:DoubleTurning} In dimension $d=1$, for $p=9$ and $\theta_\star(9)=\frac7{18}<\theta\approx0.398889$, we observe a double turning point in the $(\lambda,\mu)$ representation of the critical points. The right plot is an enlargement of the left plot near the bifurcation point.}
\end{center}
\vspace*{-0.5cm}
\end{figure}

Bifurcations taking place \emph{to the left} only if $2<p<p_0$ could sound like a paradox in view of Fig.~\ref{Fig2} (cases $\theta=0.32$ and $\theta=0.3$) as it seems that bifurcations might take place \emph{to the left} also if $7=p_0<p\lesssim 9.91109$ (see Section~\ref{Sec:critical-theta} below), but a detailed numerical analysis shows that there is in fact a double turning point in that case. An enlightening example is shown in Fig.~\ref{Fig:DoubleTurning}.

\subsection{The critical and slightly sub-critical cases for \texorpdfstring{$d=1$}{d=1}}\label{Sec:critical-theta}

For $\theta>\theta_\star$, close enough to $\theta_\star$, symmetry breaking occurs for $\lambda<d\,\theta$ (first order phase transition) if \hbox{$\kappa(p,\theta_\star)\le d\,\theta_\star=(p-2)/(2\,p)$}. Numerically, this occurs for any $p\in(2,p_\star)$ with $p_\star\approx9.91109$. 
\begin{figure}[ht]
\vspace*{-0.25cm}
\begin{center}
\includegraphics[width=5cm]{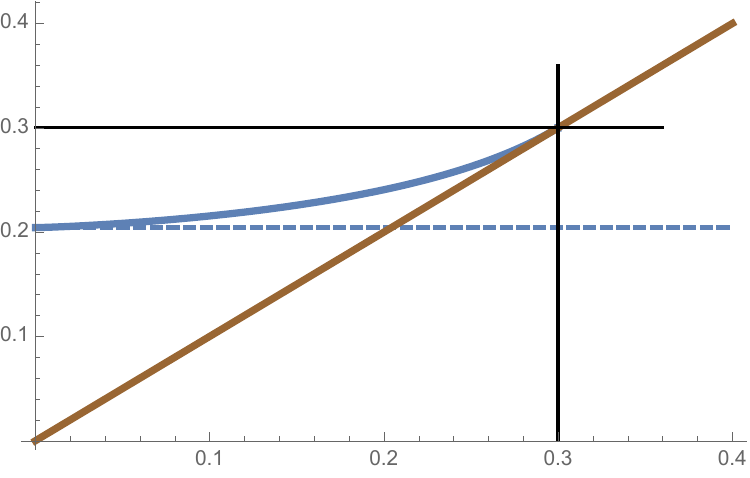}
\caption{\small\label{Fig3} In dimension $d=1$, with $p=5$ and $\theta=\theta_\star=0.3$ as in Fig.~\ref{Fig2}, we obtain numerically that $\mu(p,\theta_\star,\lambda)=\min\left\{\lambda,\kappa(p,\theta_\star)\right\}$. Proving that either $E\mapsto\lambda(\theta_\star,E)$ or $E\mapsto\mu(\theta_\star,E)$ are monotone functions are open questions. At the level $\mu=\kappa(p,\theta_\star)$, there is a loss of compactness in the variational problem associated with~\eqref{Ineq:GNS2}.}
\end{center}
\vspace*{-0.5cm}
\end{figure}

\begin{figure}[ht]
\vspace*{-0.25cm}
\begin{center}
\includegraphics[width=3.5cm]{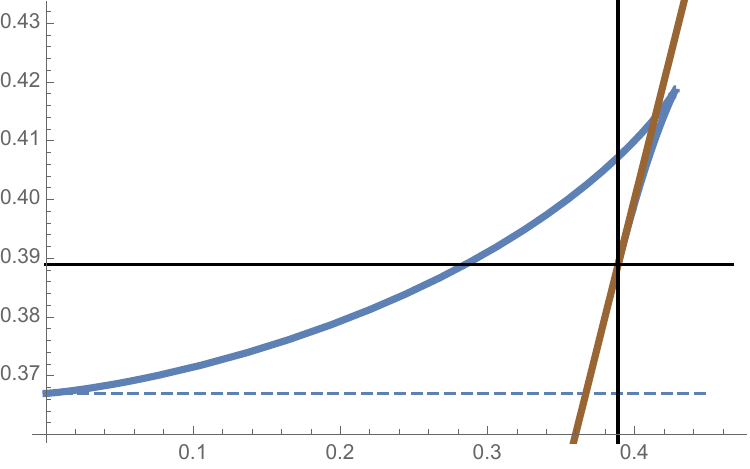}\hspace*{0.3cm}\includegraphics[width=3.5cm]{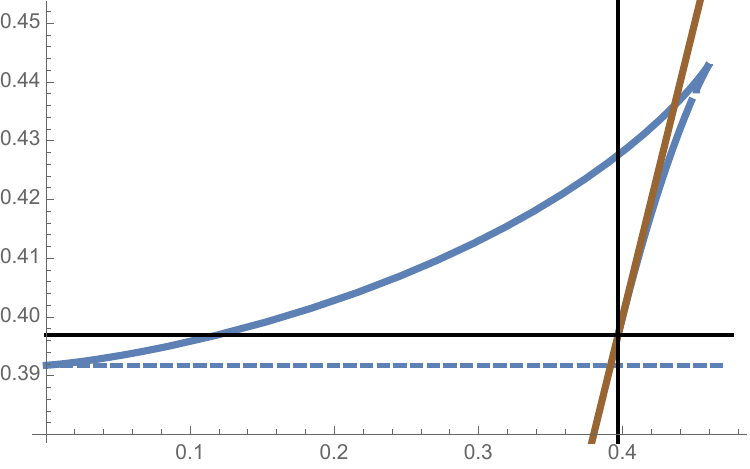}\hspace*{0.3cm}\includegraphics[width=3.5cm]{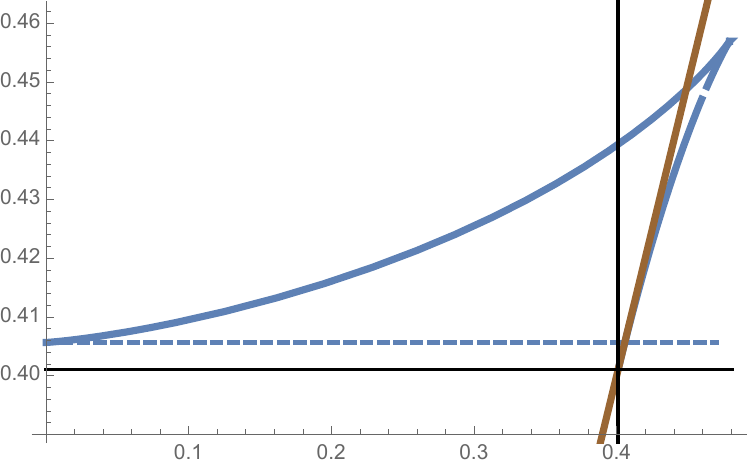}\hspace*{0.3cm}\includegraphics[width=3.5cm]{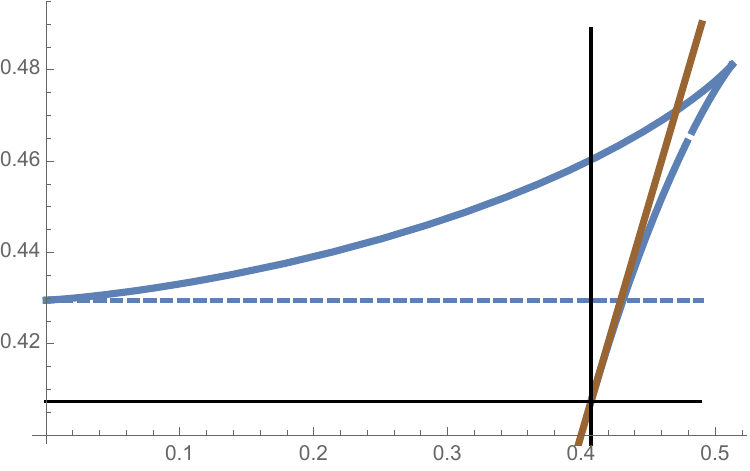}
\caption{\small\label{Fig4} In dimension $d=1$ and $\theta=\theta_\star$, the critical cases are shown (from left to right) for $p=9.0$, $9.7$, $10.1$ and $10.8$. Even in the cases of a first order transition corresponding to $\theta-\theta_\star>0$, small, and $p=9.0$ or $9.7$ the branch bifurcates \emph{to the right} at $(\lambda,\mu)=(\theta,\theta)$.}
\end{center}
\vspace*{-0.5cm}
\end{figure}

As in Section~\ref{Sec:GN}, an almost explicit range is given by taking Gaussian test functions. With the notation of~\eqref{gpd}, this amounts to the estimate $\kappa(p,\theta_\star)\le\mathsf g(p,1)=\frac{p-2}{8\,\pi}\(\frac p2\)^{2/(p-2)}$, which is not as good as the previous one, but one can check that $\mathsf g(p,1)<\theta_\star$ holds if and only if $p\in(2,p_\star)$ where $p_\star=2\,\log(2\,\pi)/\ell\approx7.8834$ where $\ell$ is the unique real number larger than $1$ such that $2\,\pi\,\ell\,e^{-\ell}=\log(2\,\pi)$. See Figs.~\ref{Fig3} and~\ref{Fig4}.

\subsection{A summary for \texorpdfstring{$d=1$}{d=1}}

\begin{figure}[hb]
\vspace*{-0.25cm}
\begin{center}
\includegraphics[width=6cm]{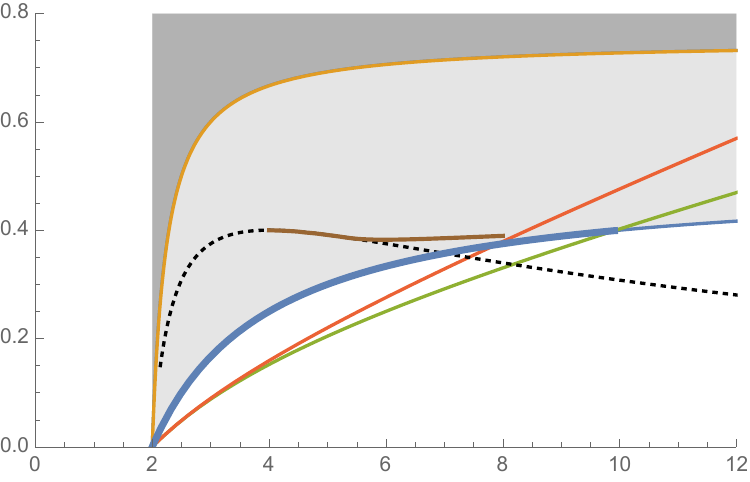}
\caption{\label{Fig5} Range of $\theta$ (vertical axis) in terms of $p$ (horizontal axis): the admissible range is the grey area $\theta\ge\theta_\star(p)$ where $p\mapsto\theta_\star(p)$ corresponds to the blue curve. The yellow curve represents $p\mapsto\theta^\#(p)$ above which the phase transition is of second order. The green curve is $p\mapsto\kappa(p,\theta_\star)$, such that symmetry breaking occurs for $\lambda<\theta$ if $\theta-\theta_\star(p)>0$ is small enough and $\kappa(p,\theta_\star)<\theta_\star(p)$, that is, in an upper neighbourhood of the bold part of the blue curve. The intersection of the red curve given by Gaussian test functions with $p\mapsto\theta_\star(p)$ determines the range $2<p<p_\ast$ in Proposition~\ref{Prop:SymmetryBreaking}, which is a more explicit interval in $p$ for symmetry breaking, for $\theta-\theta_\star>0$, small enough. The black, dotted curve is $p\mapsto\theta_0(p)$ as in Fig.~\ref{Fig:Turning}. The (brown, thick) curve $p\mapsto\theta_\bullet(p)$ is a numerical approximation of the threshold between first and second order phase transitions for a subinterval of the admissible values of $p$.}
\end{center}
\vspace*{-0.5cm}
\end{figure}

Our symmetry breaking results are summarized in~Fig.~\ref{Fig5}. A numerical estimate of the range in which a first order phase transition holds is obtained, in terms of $\theta$ by finding, for each value of $p$, for which $E\in(-\,\theta_\star,0)$ we can solve
\be{theta(E)}
\lambda(\theta,E)=\mu(\theta,E)=\theta
\ee
where $\lambda(\theta,E)$ and $\mu(\theta,E)$ are defined by~\eqref{Reparametrisation}. As in Section~\ref{Sec:Reparam}, since the solution $u_E$ of~\eqref{Eqn:lambda} corresponding to $\lambda=\overline{\lambda}(E)$ is optimal for~\eqref{Ineq:GNS1}, we know that
\[
\frac{\nrms{u_E}p^2}{\nrms{u_E}2^2}=\frac{(p-2)\,\overline{\nu}(E)+\overline{\lambda}(E)}{\overline{\mu}(E)}
\]
Altogether, we obtain
\[
\lambda(\theta,E)=\theta\,\overline{\lambda}(E)-(1-\theta)\,(p-2)\,\overline{\nu}(E)\quad\mbox{and}\quad
\mu(\theta,E)=\theta\,\big((p-2)\,\overline{\nu}(E)+\overline{\lambda}(E)\big)^{1-\frac1\theta}\,\overline{\mu}(E)^\frac1\theta\,.
\]
Solving~\eqref{theta(E)} amounts to $\theta=\widetilde\theta(E):=\frac{(p-2)\,\overline{\nu}(E)}{(p-2)\,\overline{\nu}(E)+\overline{\lambda}(E)-1}$ while $E$ is obtained as the solution of
\[
(\overline{\lambda}(E)-1)\,\log\big((p-2)\,\overline{\nu}(E)+\overline{\lambda}(E)\big)=\big((p-2)\,\overline{\nu}(E)+\overline{\lambda}(E)-1\big)\,\log\overline{\mu}(E)\,.
\]
We denote by $\widetilde E$ the numerical solution. It turns out that this problem is rather stiff and would require a more detailed numerical analysis but one can expect that $\widetilde\theta(\widetilde E)$ is the threshold value of $\theta$ below which a first order transition occurs. This value of course depends on $p$ and a portion of the curve, denoted by $p\mapsto\theta_\bullet(p)$ is shown on Fig.~\ref{Fig5}.

\subsection{Some numerical results in higher dimensions}\label{Sec:numdge3}

\begin{figure}[ht]
\vspace*{-0.5cm}
\begin{center}
\includegraphics[width=6cm]{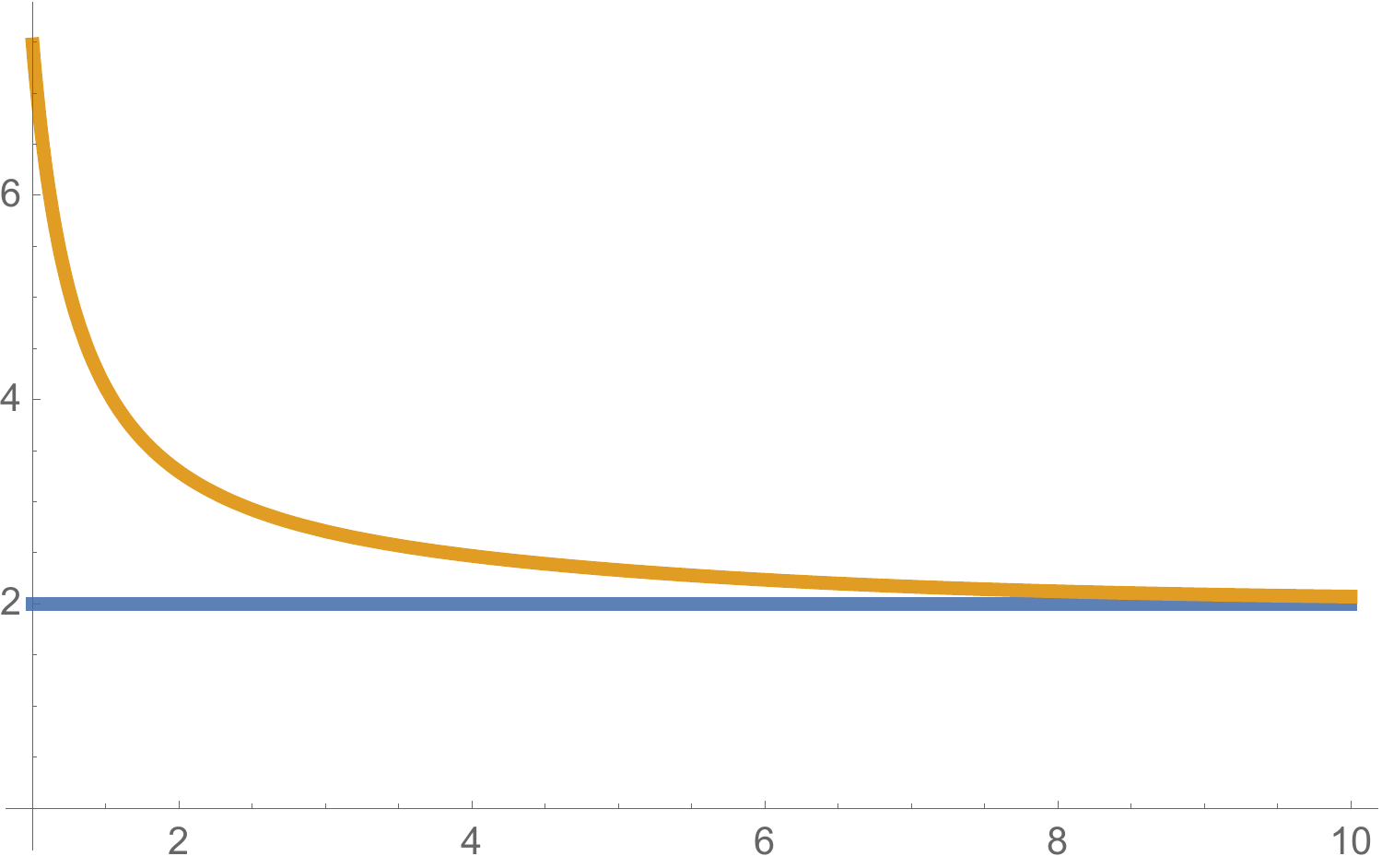}
\caption{\label{Fig:gpd} Plot (orange curve) of $d\mapsto p_\ast(d)$ where $p=p_\ast$ solves $\mathsf g(p,d)=d\,(p-2)/(2\,p)$, where $d$ is considered here as a real parameter. For integer values of $d\ge1$,  a first order phase transition occurs if $p\in(2,p_\ast(d))$.}
\end{center}
\end{figure}

The condition~\eqref{gpd} can be implemented numerically to give a range of values of $p>2$ for which  a first order phase transition occurs for some $\theta>\theta_\star$, close enough to $\theta_\star$. See Fig.~\ref{Fig:gpd}.

\newpage\noindent{\bf Acknowledgements:} This work has been partially supported by the Projects EFI (ANR-17-CE40-0030) and \emph{Conviviality} (ANR-23-CE40-0003) of the French National Research Agency. E.B.D. is funded by the European Union’s Horizon 2020 research and innovation programme under the Marie Skłodowska-Curie grant agreement n\textsuperscript{o} 101034255.\\[-2pt]
\noindent{\scriptsize \copyright\,\the\year~by the authors. Reproduction of this article by any means permitted for non-commercial purposes. \hbox{\href{https://creativecommons.org/licenses/by/4.0/legalcode}{CC-BY 4.0}}}
\bibliographystyle{siam}\small\begin{spacing}{1.05}
\bibliography{SBI}
\end{spacing}
\end{document}